\documentclass[11pt,reqno]{amsart}

\usepackage[text={160mm,240mm},centering]{geometry}            % See geometry.pdf to learn the layout options. There are lots.
\usepackage{amssymb,amsmath,setspace,geometry,indentfirst,changepage,inputenc,amsthm}
\usepackage{mathrsfs,amsfonts,savesym,graphicx,bm,color}

\usepackage{graphicx}
\usepackage{amssymb}
\usepackage{epstopdf}
\usepackage{dsfont}

\usepackage{lscape}

\usepackage{tikz}

\usepackage{float}

\usepackage[all,2cell]{xy}    % ???? xy ø?
\UseAllTwocells               % ????????? 2-cell

\usepackage{amsmath,amsthm}

\usepackage[mathscr]{eucal}
\usepackage[all]{xy}
\usepackage{mathrsfs}
\usepackage{hyperref}
\usepackage{color}

\newtheorem{theorem}{Theorem}[section]
\newtheorem{lemma}[theorem]{Lemma}

\newtheorem{thmletter}{Theorem}

\newtheorem{conjecture}[theorem]{Conjecture}
\newtheorem{corollary}[theorem]{Corollary}
\newtheorem{proposition}[theorem]{Proposition}

\newtheorem{fact}[theorem]{Fact}
\newtheorem{question}[theorem]{Question}
\newtheorem{definition-lemma}[theorem]{Definition-Lemma}
\newtheorem{definition-theorem}[theorem]{Definition-Theorem}

\newtheorem{algorithm}[theorem]{Algorithm}
\newtheorem{problem}[theorem]{Problem}

\theoremstyle{definition}
\newtheorem{example}[theorem]{Example}
\newtheorem{definition}[theorem]{Definition}

\newtheorem{remark}[theorem]{Remark}

\allowdisplaybreaks

\title[]{On jet closures of singularities}

\author{Yifan Chen}
\address{
	Department of Mathematical Sciences,
	Tsinghua University,
	Beijing, 100084, P. R. China.}
\email{c-yf20@tsinghua.org.cn}
\author{Huaiqing Zuo}
\address{Department of Mathematical Sciences,
	Tsinghua University,
	Beijing, 100084, P. R. China.}
\email{hqzuo@mail.tsinghua.edu.cn}

\begin{document}
	
	\maketitle
	
	\parskip=.3em
	
	\setcounter{tocdepth}{1}
	%\tableofcontents

\begin{abstract}
	The jet closure and jet support closure were first introduced by de Fernex, Ein and Ishii to solve the local isomorphism problem. In this paper, we introduce two local algebras associated to jet closure and jet support closure respectively. We show that these two algebras are invariants of the singularities. We compute and investigate these invariants for some interesting cases, such as the cases of monomial ideals and homogeneous ideals. We also introduce a new filtration and jet index to jet closures. The jet index describes which jet scheme recover the information of base scheme. Moreover, we obtain some properties of the jet index.

\end{abstract}

\section{Introduction}
\label{s1}

In 1968, Nash introduced arc spaces and jet schemes to singularities in \cite{7}. Later the studies of jet schemes were used to motivic integration by Denef and Loeser in \cite{8} and  Mustata used these results to study the singularities of pairs in \cite{9}. Ein and Mustata also utilized arc spaces and jet schemes to study other invariants of singularities, some results were collected in \cite{10}.  One may also refer to the articles \cite{5.1,  5.2, 5.3, 5.4} for some recent progress about singularities and jet schemes.

Let $X$ be a scheme and denote $X_{m}$ as its $m$-th jet scheme, we have a canonical projection $\pi_{m}: X_{m} \longrightarrow X$. For every $x \in X$ we denote by $X_{m}^{x}:=\pi_{m}^{-1}(x)$ the fiber over $x$ and we say it is the $m$-th scheme of local jets. In \cite{1}, Tommaso de Fernex, Lawrence Ein and Shihoko Ishii introduced jet closure of an ideal $I$ in a local ring $R$. Generally speaking, the $m$-th jet closure of an ideal $I$ in a local ring $R$, denoted by $I^{m-jc}$, is the largest ideal $J$ such that the $m$-th scheme of local jets of $Spec(R)$ defined by $J$ is equal to that defined by $I$. It was shown in \cite{1} that the jet closure is a closure operation and most time it is a non-trivial operation. When the $m$-th jet closure of an ideal $I$ is equal to $I$ itself, we say the ideal is $m$-th jet closed. The intersection of all $m$-th jet closure of an ideal $I$ is defined to be the arc closure of $I$, and similarly we say an ideal is arc closed if the arc closure of it is equal to itself. The authors of \cite{1} also introduced another closure operation, namely jet support closure $I^{m-jsc}$, which is similar to jet closure and is connected with the integral closure of an ideal with the following proposition.

\begin{proposition}
(\cite[Proposition 4.12]{1}) 
Suppose $R$ is a local integral domain essentially of finite type over a field $k$ and $I \subset R$ is an ideal. Denote $I^{m-jsc}$ the $m$-th jet support closure of $I$, and $$I^{\infty -jsc}:=\displaystyle\bigcap_{m \in \mathbb{N}}I^{m-jsc}$$ the arc support closure of $I$, then we have: 

(1)$I^{\infty -jsc} \subset \bar{I}$,  

(2)If $R$ is regular, then $I^{\infty -jsc}=\bar{I}$, 

where $\bar{I}$ denote the integral closure of $I$.
\end{proposition}

The authors also raised questions about when an ideal is arc closed in \cite{1}, and a positive answer about the ideals in Notherian rings was given in \cite{3}.

\begin{proposition}
(\cite[Theorem 5.1]{3}) 
Let $(R,\mathbf{m})$ be a local Noetherian $k$-algebra with residue field $L$. Suppose $L$ is separable over $k$, then $I=I^{ac}:=\bigcap_{m \in \mathbb{N}}I^{m-jc}$, where $I^{ac}$ denotes the arc closure of $I$.
\end{proposition}

In \cite{1}, the jet closure and jet support closure are used to characterize the local isomorphism problem, but they are of good value to study themselves because intuitively they characterize how much information is losing while using the method of jet schemes. Both the authors of \cite{1} and \cite{3} are focusing on the generic properties of jet closure, and they study mainly on arc closure. In this paper, however, we give a new point of view by introducing the algebra associated to jet closure and jet support closure which are invariants of singularities. In appendix, we will provide a computer procedure in Singular to compute these invariants. Through computing concrete examples, we find some interesting results. Some conjectures about jet closure and jet support closure are raised, and are proved strictly. 

Studying the filtration properties of jet closure is another brand-new way. We will give several monotony properties for jet closures, proving that the jet closures can define a filtration as in \cite{6}. By the monotony properties we can prove that for any Artinian algebra $k[[x_{1},...,x_{n}]]/I$ there exists an integer $s$ called "jet index" such that $I$ is $s$-jet closed. It characterizes which jet scheme recovers the information of the original scheme and it is a good new object to study with. 

Firstly we obtain some explicit results for jet closures under the homogeneous case.  In the following, throughout the paper, Let k be an algebraically closed field of characteristic 0.

\begin{thmletter}\label{mtha}
Let $R=k[[x_{1},...,x_{n}]]$ (i.e.,  a formal  power series ring)) and $I \subset R$ an ideal. If $I$ is an ideal extended from a homogeneous ideal in $k[x_{1},...,x_{n}]$, i.e. $I=(f_{1},...,f_{s})$, where $f_{1},...,f_{s}$ are homogeneous polynomials from $k[x_{1},...,x_{n}]$, then the $m$-th jet closure of $I$ is $I^{m-jc}=I+\mathbf{m}^{m+1}$, where $\mathbf{m}$ denotes the maximal ideal in $R$.
\end{thmletter}

In the practice of calculating, we find that in the local ring $k[[x,y]]$, when $f$ is a weighted homogeneous polynomial, the jet closures of $J(f)$ always have the same form, so we propose the following conjecture.
\begin{conjecture}
If $f \in k[[x_{1},...,x_{n}]]$ is a weighted homogeneous polynomial, then we have $J(f)^{m-jc}=J(f)+\mathbf{m}^{m+1}$, where $\mathbf{m}$ is the maximal ideal of $k[[x_{1},...,x_{n}]]$.
\end{conjecture}

Recall that the most important classification result for hypersurface singularities in characteristic zero is the following result by Arnold:
\begin{theorem}(Arnold; \cite{AVZ})
 Let $f \in k[[x_{1},...,x_{n}]]$  (i.e., a formal power series ring) be a simple singularity. Then $f$ is contact equivalent to an $ADE$ singularity from the following list:
 
$ \begin{array}{lllll}
A_{m}: & x_1^{m+1}+x_2^2 & +q, & m \geq 1, & q:=x_3^2+\cdots+x_n^2 \\
D_{m}: & x_2\left(x_1^2+x_2^{m-2}\right) & +q, & m \geq 4 & \\
E_6: & x_1^3+x_2^4 & +q & & \\
E_7: & x_1\left(x_1^2+x_2^3\right) & +q & & \\
E_8: & x_1^3+x_2^5 & +q & &
\end{array}
$
\end{theorem}

Arnold's classification has numerous applications. The list of simple or $ADE$ singularities appears in many other contexts of mathematics and is obtained also by classifications using a completely different equivalence relation (cf. \cite{Du}). In \cite[Proposition 5.1]{1}, the authors show that if the ideal is arc closed, one can use all of its jet closures to distinguish different germs. This criterion meets some difficulties in explicit calculation because it involves infinitely many jet closures. This paper, however, tries to distinguish different germs by finitely many jet closures or jet support closures when it is simple singularity. In general, the calculation of jet closures is not easy, and we can get the following results of simple curve singularities when $m$ is small.

\begin{thmletter}\label{mthb}
Let $R=k[[x,y]]$, $I \subset R$ defines a simple curve singularity.

(1)$A_{n-1}$, $I=(x^{2}+y^{n})$ $(n \ge 3)$. When $m \le n$, $I^{m-jc}=I+\mathbf{m}^{m+1}$, and when $m=n+1$, $I^{m-jc}=\mathbf{m}^{m+1}+I+(x^{3})$.  

(2)$D_{n}$, $I=(x^{2}y+y^{n-1})$. When $m \le n-1$, $I^{m-jc}=I+\mathbf{m}^{m+1}$.  

(3)$E_{6}$, $I=(x^{3}+y^{4})$. When $m \le 4$, $I^{m-jc}=I+\mathbf{m}^{m+1}$.

(4)$E_{7}$, $I=(x^{3}+xy^{3})$. When $m \le 4$, $I^{m-jc}=I+\mathbf{m}^{m+1}$.

(5)$E_{8}$, $I=(x^{3}+y^{5})$. When $m \le 5$, $I^{m-jc}=I+\mathbf{m}^{m+1}$. 
\end{thmletter}

For jet support closures, the situation becomes more complicated. However, when the ideal is a monomial ideal, the jet support closure can also be calculated explicitly, though in a more complicated way.

\begin{thmletter}\label{mthc}
Let $R=k[[x_{1},...,x_{n}]]$, for a monomial ideal $I \subset R$, the $m$-th jet support closure of I is also a monomial ideal, and it consists of the following monomials: $x_{1}^{b_{1}}...x_{n}^{b_{n}} \in I^{m-jsc}$ if and only if for any $t_{1},...,t_{n} \in \mathbb{N}$ satisfying $t_{1}b_{1}+t_{2}b_{2}+...+t_{n}b_{n} \le m$, there exists $x_{i_{1}}^{a_{1}}...x_{i_{s}}^{a_{s}} \in I$ such that $t_{i_{1}}a_{1}+...+t_{i_{s}}a_{s} \le m$. 
\end{thmletter}

When the ideal is generated by a single polynomial $f$, the jet support closure of $(f)$ can be computed in some special cases. First we give the definition of weighted homogeneous polynomial and the corresponding weight.

\begin{definition}
A polynomial $f(x_{1},...,x_{n})$ is called a weighted homogeneous polynomial if there exist $d_{1},...,d_{n} > 0$ such that for every monomial $x_{1}^{a_{1}}...x_{n}^{a_{n}}$ in $f$, $a_{1}d_{1}+...+a_{n}d_{n}=d$ for some fixed $d$. The numbers $d_{1},...,d_{n}$ are called the weight of $x_{1},...,x_{n}$.
\end{definition}

We calculate all the jet support closure for weighted homogeeous polynomials in two variables. In application, we have the jet support closure for simple curve singularities, and we can distinguish different simple curve singularities by finitely many jet support closures. The number of jet support closures we need is near the Milnor number of the singularity, which we give the definition below:
\begin{definition}
If $f \in k[[x_{1},...,x_{n}]]$  defines an isolated hypersurface singularity, then we call the $k$-dimension of $k[[x_{1},...,x_{n}]]/(\partial_{x_{1}}f,...,\partial_{x_{n}}f)$ the Milnor number of $k[[x_{1},...,x_{n}]]/(f)$, denoted by $\mu(k[[x_{1},...,x_{n}]]/(f))$. 
\end{definition}

\begin{thmletter}\label{mthd}
(1)Let $f(x_{1},...,x_{n})$ be a homogeneous polynomial of degree $d$, suppose $f$ is reduced, then $(f)^{m-jsc}=(f)+\mathbf{m}^{m+1}$ in $k[[x_{1},...,x_{n}]]$. 

(2)Suppose $f \in R=k[[x,y]]$ is a weighted homogeneous polynomial, the weight of $x$ is $a$ and the weight of $y$ is $b$, satisfying $a,b \in \mathbb{N}$ and $gcd(a,b)=1$. We assume $f(x,y)=\sum_{l=1}^{s}c_{l}x^{i_{l}}y^{j_{l}}$, where $ai_{l}+bj_{l}=d$ and $c_{l} \in k/\{0\} $ for every $l$. Let $A_{m}:=\{ (u,v) \in \mathbb{N}^{2}: ui_{l}+vj_{l}\ge m+1, l=1,...,s \}$, then when $m<d$, $(f)^{m-jsc}$ is generated by all $x^{p}y^{q}$, where $(p,q)$ satisfies $pu+qv \ge m+1$ for all $(u,v) \in A_{m}$, and when $m \ge d$, $I^{m-jsc}=(f,x^{p_{1}}y^{q_{1}})\cap(x^{p_{2}}y^{q_{2}})$, where $ap_{1}+bq_{1} \ge m+1$ and $up_{2}+vq_{2} \ge m+1$ for all $(u,v) \in A_{m}, u \le a \ or \ v \le b$. 

(3)Suppose $R=k[[x,y]]$, and $R/I_{1}$ and $R/I_{2}$ are two simple curve singularities. Set $M:=max\{\mu(R/I_{1}),\mu(R/I_{2})\}$ is the maximum of  Milnor numbers. Then we have $R/I_{1} \cong R/I_{2}$ if and only if $R/I_{1}^{m-jsc} \cong R/I_{2}^{m-jsc}$ for all $m \le M+1$.
\end{thmletter}

This paper is organized as follows. In section \ref{s2}, we review basic properties of the jet scheme, jet closure and jet support closure.  We will demonstrate how to compute them explicitly. We will introduce an algebra associated to a  jet closure. We show that this algebra is a invariant of singularities in section \ref{s3}, moreover, we will study them in several interesting cases, proving Theorem \ref{mtha} and Theorem \ref{mthb}. In section \ref{s4}, we will present the properties  of jet support closures of monomial ideals, by proving Theorem \ref{mthc}. In section \ref{s5}, the jet support closure of a reduced homogeneous polynomial is computed. We give the proof of Theorem \ref{mthd}. In section \ref{s6}, we  investigate the jet closure as a filtration. Moreover, we define the jet index and raise some questions about it. At last, we provide two singular programs for the computing jet closure and jet support closure respectively in section \ref{s7}.

\section{Preliminaries}
\label{s2}

\subsection{Jet scheme}
\label{s2.1}

In this section, we will introduce some basics about jet schemes and use Hasse-Schmidt derivation to calculate them. More comprehensive introduction can be found in, for example, \cite{4} and \cite{5}. First we will give the definition of jet scheme.

\begin{definition}
Let $X$ be a scheme over a field $k$, for every $m \in \mathbb{N}$, consider the functor from $k$-schemes to set
\begin{flalign*}
Z \mapsto Hom(Z \times_{Spec(k)} Spec (k[t]/(t^{m+1})),X).
\end{flalign*}
There is a k-scheme representing this functor (see for
example \cite[Theorem 2.1]{5}), which is called the $m$-th jet scheme of $X$, denoted by $X_{m}$, i.e.
\begin{flalign*}
Hom(Z,X_{m})=Hom(Z \times_{Spec(k)} Spec (k[t]/(t^{m+1})),X).
\end{flalign*}
\end{definition}

For $1 \le i \le j$, the truncation map $k[t]/(t^{j}) \longrightarrow k[t]/(t^{i})$ induces natural projections between jet schemes $\psi_{i,j}: X_{j} \longrightarrow X_{i}$. If we identify $X_{0}$ with $X$, we can get natural projection $\pi_{i}: X_{i} \longrightarrow X$. One can check that $\{X_{m}\}$ is an inverse system and the inverse limit $X_{\infty}:=\lim \limits_{\leftarrow m}X_{m}$ is called the arc space of $X$. There are also natural projections $\psi_{i}: X_{\infty} \longrightarrow X_{i}$. If $X$ is of finite type, $X_{m}$ is also of finite type for each $m \in \mathbb{N}$, but $X_{\infty}$ is usually not. Jet scheme is also compatible with homomorphisms between $k$-schemes, with the following proposition.

\begin{proposition}
(\cite[Proposition 2.4]{5}) 
Let $f: X \longrightarrow Y$ be a morphism of $k$-schemes of finite type, then there is a natural morphism $f_{m}: X_{m} \longrightarrow Y_{m}$ such that the following diagram is commutative: \\
\begin{center}
\
\xymatrix@=8ex
{
X_{m} \ar[r]^{f_{m}} \ar[d]^{\pi_{X_{m}}} & Y_{m} \ar[d]^{\pi_{Y_{m}}} \\
X \ar[r]^{f} & Y.
} \\
\end{center}
\end{proposition}

In order to compute jet scheme explicitly, we introduce the Hasse-Schmidt derivation.
\begin{definition}
Let $B$,$R$ be $k$-algebras, and $m \in \mathbb{N}$, a Hasse-Schmidt derivation of order $m$ from $B$ to $R$ over $k$ is a sequence $(D_{0},D_{1},...,D_{m})$ such that $D_{0}$ is a homomorphism of $k$-algebras, and each $D_{i}(1 \le i \le m)$ is a homomorphism of abelian groups such that: 

(1)$D_{i}(a)=0$ for each $a \in k, 1 \le i \le m$,

(2)For all $x,y \in B$, $0 \le k \le m$,
\begin{flalign*}
D_{k}(xy)=\displaystyle\sum_{i+j=k}D_{i}(x)D_{j}(y).
\end{flalign*}

The set of Hasse-Schmidt derivation of order $m$ from $B$ to $R$ over $k$ is denoted by $Der_{k}^{m}(B,R)$ or $Der^{m}(B,R)$ for short.
\end{definition}

The following definition gives an explicit construction of jet scheme (we will give the proposition later).

\begin{definition}
Let $B$ be a $k$-algebra and $m \in \mathbb{N}$, define a $B$-algebra $HS_{B}^{m} :=B[x^{(i)}]_{x \in B, i=1,2,...,m}/I$, where $I$ is generated by the following set: 

(1)$\{ (x+y)^{(i)}-x^{(i)}-y^{(i)}:x,y \in B,i=1,2,...,m\}$. 

(2)$\{a^{(i)}:a \in k,i=1,...,m \}$.

(3)$\{(xy)^{(k)}-\displaystyle\sum_{i+j=k}x^{(i)}y^{(j)}:x,y \in B, k=0,...,m\}$,
where we identify $x^{(0)}$ with $x$ for all $x \in B$. 

We also define the universal derivation $(d_{0},...,d_{m})$ from $B$ to $HS_{B}^{m}$ by $d_{i}(x)=x^{(i)}$.
\end{definition}

The following two propositions connect the $Der^{m}(B,R)$ and $HS_{B}^{m}$ with the $k$-morphisms between $B$ and $R[t]/(t^{m+1})$. They are proven in \cite{4}.

\begin{proposition}
(\cite[Lemma 1.7]{4}) 
For any $(D_{0},...,D_{m}) \in Der^{m}(B,R)$, define $\phi :B \longrightarrow R[t]/(t^{m+1})$ such that $x \mapsto D_{0}(x)+...+D_{m}(x)t^{m}$, then we have $\phi \in Hom(B,R[t]/(t^{m+1}))$ and the correspondence $Der^{m}(B,R) \longrightarrow Hom(B,R[t]/(t^{m+1}))$ is bijective.
\end{proposition}

\begin{proposition}(\cite[Corollary 1.8]{4}) 
For any $\phi \in Hom(HS_{B}^{m},R)$, define a morphism from $B$ to $R[t]/(t^{m+1})$ given by $x \mapsto \phi(d_{0}(x))+...+\phi(d_{m}(x))t^{m}$, and this gives a bijection between $Hom(HS_{B}^{m},R)$ and $Hom(B,R[t]/(t^{m+1}))$.
\end{proposition}

The proposition above shows that when the scheme $X=Spec(B)$, $Spec(HS_{B}^{m})$ is the required $m$-th jet scheme, i.e. $X_{m}=Spec(HS_{B}^{m})$. For the convenience of notation, when the scheme $X=Spec(B)$, we will denote $HS_{B}^{m}$ by $B_{m}$, $X_{m}=Spec(B_{m})$ and we say $B_{m}$ is the $m$-th jet ring. In particular, when $B=R/I$, where $I \subset R$ is an ideal, we have $B_{m}=R_{m}/I_{m}$, and $I_{m}$ is called $m$-th jet ideal. Here is an example of explicit computation.

\begin{example}
Let $X=Spec(R)$ and $R=k[x,y]/(x^{2}-y^{3})$. To compute the $2$-th jet scheme of $X$, let $x=x_{0}+x_{1}t+x_{2}t^{2}$ and $y=y_{0}+y_{1}t+y_{2}t^{2}$. Set $f(x,y)=x^{2}-y^{3}$, and calculate $f(x_{0}+x_{1}t+x_{2}t^{2},y_{0}+y_{1}t+y_{2}t^{2})$, we can get $R_{2}=k[x_{0},x_{1},x_{2},y_{0},y_{1},y_{2}]/I_{2}$, where $I_{2}=(x_{0}^{2}-y_{0}^{3},2x_{0}x_{1}-3y_{0}^{2}y_{1},2x_{0}x_{2}+x_{1}^{2}-3y_{0}^{2}y_{2}-3y_{1}^{2}y_{0})$. We then get the $2$-th jet ring of $R$ and the $2$-th jet ideal of $I$, namely $R_{2}$ and $I_{2}$.
\end{example}

\subsection{Jet Closure and Jet Support Closure}
\label{s2.2}

In 2017, Tommaso de Fernex, Lawrence Ein, Shihoko Ishii introduced jet closure and jet support closure in \cite{1}. In this section, we will recall some basic properties of them. Firstly, we fix a local ring $(R, \mathbf{m},k)$, let $R_{m}$ denotes the $m$-th jet ring and $I_{m}$ denotes the m-th jet ideal for any ideal $I \subset R$. Let $X$=$Spec(R)$, and $\mathit{o} \in X$ is the closed point. Assume $X_{m}=Spec(R_{m})$ is the $m$-th jet scheme of $X$. If $V(I) \subset X$ is the subscheme defined by $I$, then its $m$-th jet scheme $V(I)_{m}$ is the subscheme defined by $I_{m}$ in $X_{m}$. Let $\varphi : X_{m} \longrightarrow X$ be the natural projection, and $X_{m}^{o}$ denotes the fiber over $\mathit{o}$ and $V(I)_{m}^{o}$ is the fiber of $\mathit{o}$ in the projection $V(I)_{m} \longrightarrow V(I)$.

\begin{definition}
For any ideal $I\subset R$ and $m\in \mathbb{N} \cup \{ \infty \} $, we define the $m$-th jet closure of $I$ to be
\begin{flalign*}
I^{m-jc}:= \{ f \in R : (f)_{m} \subset I_{m} (mod  \ \mathbf{m}R_{m}) \}.
\end{flalign*}
For $m=\infty$, we also call the ideal $I^{\infty -jc}$ the arc closure of $I$.
\end{definition}

Here are some basic properties of jet closure which are proved in \cite{1}.

\begin{fact}
(\cite[Lemma 3.2]{1})
Assume $\phi : R \longrightarrow R/I$ is the quotient map, then $I^{m-jc}$ in R equals to $(0)^{m-jc}$ in $R/I$.
\end{fact}

\begin{fact}
(\cite[Lemma 3.3]{1})
We define a ring map $\gamma_{m}:R \longrightarrow R_{m}/ \mathbf{m} R_{m}[t]/(t^{m+1})$ such that $\gamma_{m}(f)= \sum_{i=0}^{m}d_{i}(f)t^{i}$ where $d_{i}(f)$ denotes the universal Hasse-Schmidt derivations of the element of $f$ in $I$, then $(0)^{m-jc}=ker(\gamma_{m})$.
\end{fact}

According to the two facts above, we have the following corollary.

\begin{corollary}
Suppose $I \subset R$ is an ideal ,and if we define the ring map to be $\mu_{m}:R \longrightarrow R_{m}/ \mathbf{m} R_{m}[t]/(t^{m+1},I_{m})$ where $\mu_{m}(f)= \sum_{i=0}^{m}d_{i}(f)t^{i}$, then $(I)^{m-jc}=ker(\mu_{m})$.
\end{corollary}

Using the corollary above, we can compute the jet closure explicitly through Singular when $R=k[[x_1,...,x_{n}]]$ and the code can be found in the appendix. Here we present an example to show how to compute jet closure by hand, which is useful in the proof in the next several sections.

\begin{example}
We continue the Example $2.1$, but this time we set $R=k[x,y]$ and $I=(x^{2}-y^{3})$. We have computed $I_{2}=(x_{0}^{2}-y_{0}^{3},2x_{0}x_{1}-3y_{0}^{2}y_{1},2x_{0}x_{2}+x_{1}^{2}-3y_{0}^{2}y_{2}-3y_{1}^{2}y_{0})$. To see if an arbitrary $f \in R$ is in $I^{2-jc}$, it suffices to see if $\mu_{2}(f)$ is in $(t^{3},I_{2},\mathbf{m})$. Thus we need to compute $f(x_{1}t+x_{2}t^{2},y_{1}t+y_{2}t^{2})=f_{1}t+f_{2}t^{2}+...$ and $f \in I^{2-jc}$ if and only if $f_{1},f_{2} \in I_{2}$. Here when we compute $\mu_{2}(f)$ we "omit" the constant term because they are already in the $(t^{3},I_{2},\mathbf{m})$.
\end{example}

The next fact shows the jet closure is indeed a closure operation.

\begin{fact}
(\cite[Corollary 3.6]{1})
For any ideal $I \subset R$ and any $m \in \mathbb{N}$ we have $(I^{m-jc})^{m-jc}=I^{m-jc}$.
\end{fact}

\begin{definition}
We say an ideal $I$ is $m$-jet closed if $I=I^{m-jc}$, and $I$ is arc closed if $I^{\infty -jc}=I$.
\end{definition}

The following fact shows that in general, jet closure is a nontrivial operation.

\begin{fact}
(\cite[Proposition 3.9]{1})
For any ideal $I \subset R$, we have $I+\mathbf{m}^{m+1} \subset I^{m-jc}$.
\end{fact}

The next property characterizes the geometric meaning of $m$-jet closures.

\begin{fact}
(\cite[Proposition 3.4]{1})
For any ideal $I \subset R$, $I^{m-jc}$ is the largest ideal $J \subset R$ such that $V(I)_{m}^{o}=V(J)_{m}^{o}$.
\end{fact}

The basic connection between jet closure and arc closure is the following fact.

\begin{fact}
(\cite[Proposition 3.12]{1})
For any ideal $I \subset R$, $\bigcap_{m \in \mathbb{N}}I^{m-jc}=I^{\infty -jc}$.
\end{fact}

In \cite{1}, the jet closure and arc closure are used to solve the embedded local isomorphism problem. Here is the setting. Let $\phi:(Y,y) \longrightarrow (X,x)$ be a $k$-morphism of $k$-schemes, then there are  induced morphisms: $\phi_{m}:Y_{m} \longrightarrow X_{m}$ and the local morphisms $\phi_{m}^{loc}:Y_{m}^{y} \longrightarrow X_{m}^{x}$, where $X_{m}^{x}$ and $Y_{m}^{y}$ denote the fiber of $X$ and $Y$ over $x$ and $y$. There are two problems called local isomorphism problem and embedded local isomorphism problem. 

\begin{problem}
(Local isomorphism problem)  
With the notations above, if $\phi_{m}^{loc}$ is an isomorphism for every $m \in \mathbb{N} \cup \{ \infty \}$, does it follow that $\phi$ is an isomorphism of germs?
\end{problem}

\begin{problem}
(Embedded local isomorphism problem) 
With the notations above, if $\phi$ is a closed immersion and $\phi_{m}^{loc}$ is an isomorphism for every $m \in \mathbb{N} \cup \{ \infty \}$, does it follow that $\phi$ is an isomorphism of germs?
\end{problem}

\begin{definition}
Given a germ $(X,x)$. If Problem $2.1$(resp.Problem $2.2$) has a positive answer for a certain class of germs $(Y,y)$, then we say $(X,x)$ has the local isomorphism property(resp.embedded local isomorphism property) for that certain class of germs.
\end{definition}

In \cite{1}, there is a proposition which connects the two problems.

\begin{proposition}
(\cite[Proposition 2.6]{1}) 
Suppose $(X,x)$ has the embedded local isomorphism property for every germs $(Y,y)$, then $(X,x)$ has the local isomorphism property for all locally Noetherian $k$-scheme germs $(Y,y)$.
\end{proposition}

The equivalence of the study of arc closure and the embedded local isomorphism property was shown in \cite{1}.

\begin{theorem}
(\cite[Proposition 5.1]{1}) 
Assume $(R, \mathbf{m})$ is a local $k$-algebra. Let $X=Spec(R)$ and $\mathit{o} \in X$ is the closed point, then $(X,\mathit{o})$ has the embedded isomorphism property for every $(Y,y)$ if and only if the zero ideal in $R$ is arc closed.
\end{theorem}

In \cite{1}, the authors gave a question about the arc closure: when is an ideal arc closed? The question was answered in \cite{3}, which tells us in general the answer is positive.

\begin{theorem}
(\cite[Theorem 5.1]{3}) 
Let $(R, \mathbf{m})$ be a Noetherian local $k$-algebra with residue field $L$, if $L/k$ is separable, then for any proper ideal $I \subset R$, $I$ is arc closed.
\end{theorem}

Another concept similar to jet closure is jet support closure, with the following definition.

\begin{definition}
For any ideal $I\subset R$ and $m\in \mathbb{N} \cup \{ \infty \} $, we define the $m$-th jet support closure of $I$ to be
\begin{flalign*}
I^{m-jsc}:= \{ f \in R | (f)_{m} \subset \sqrt{I_{m}+\mathbf{m}R_{m}} \}.
\end{flalign*}
For $m=\infty$, we also call the ideal $I^{\infty -jsc}$ the arc support closure of $I$.
\end{definition}

The jet support closure has many properties similar to which of jet closure, and we summarize them as the following proposition.

\begin{proposition}
(Basic properties of jet support closure, \cite[Proposition 4.3, 4.4, 4.7]{1}) 

(1)Assume $\phi : R \longrightarrow R/I$ is the quotient map, then $I^{m-jsc}$ in $R$ equals to $(0)^{m-jsc}$ in $R/I$. 

(2)We define a ring map $\nu_{m}:R \longrightarrow (R_{m}/ \mathbf{m} R_{m})_{red}[t]/(t^{m+1})$ such that $\nu_{m}(f)= \displaystyle\sum_{i=0}^{m}d_{i}(f)t^{i}$ where $d_{i}(f)$ denotes the universal Hasse-Schmidt derivations of the element of $f$ in $I$, then $(0)^{m-jsc}=ker(\nu_{m})$. 

(3)Suppose $I \subset R$ is an ideal ,and if we define the ring map to be $\lambda_{m}:R \longrightarrow R_{m}[t]/(t^{m+1},\sqrt{I_{m}+\mathbf{m}^{e}})$ where $\lambda_{m}(f)= \displaystyle\sum_{i=0}^{m}d_{i}(f)t^{i}$, then $(I)^{m-jsc}=ker(\lambda_{m})$. 

(4)$(I^{m-jsc})^{m-jsc}=I^{m-jsc}$. 

(5)$I^{m-jc} \subset I^{m-jsc}$. 

(6)$I^{m-jsc}$ is the largest ideal $J \subset R$ such that $(V(I)_{m}^{o})_{red}=(V(J)_{m}^{o})_{red}$. 
\end{proposition}

\section{Calculation of Specific Jet Closures}
\label{s3}

In this section, we will fix the local ring $R=k[[x_1,...,x_{n}]]$ and we will introduce some properties for the computation of jet closures. First we will prove that the jet closure is an invariant of singularity.

\begin{proposition}
Fix $R=k[[x_1,...x_{n}]]$, $I,J \subset R$ are ideals, assume there exists an isomorphism $\phi:R \longrightarrow R$ such that $\phi(I)=J$, then $\phi$ induces $R/I^{m-jc} \cong R/J^{m-jc}$ for any integer $m$, $m \ge 1$. 
\end{proposition}
\begin{proof}
Let $I_{m}$ be the m-th jet ideal of $I$ contained in $R_{m}=k[[x_{i}^{(j)}]](1 \le i \le n,0 \le j \le m)$, $\mathbf{m} \subset R$ be its maximal ideal and $\mathbf{m}^{e}$ be its extension in $R_{m}$, then according to Corollary $2.1$, there exists a ring map $\mu_{m}:R \longrightarrow R_{m}[t]/(t^{m+1},I_{m},\mathbf{m}^{e})$,  
where $\mu_{m} (x_{i})=x_{i}^{(0)}+x_{i}^{(1)}t+...+x_{i}^{(m)}t^{m}$, and  $ker(\mu_{m})=I^{m-jc}$. Similarly, there exists a ring map $\nu_{m}:R \longrightarrow R_{m}[t]/(t^{m+1},J_{m},\mathbf{m}^{e})$, 
where $\nu_{m} (x_{i})=x_{i}^{(0)}+x_{i}^{(1)}t+...+x_{i}^{(m)}t^{m}$, and  $ker(\nu_{m})=J^{m-jc}$. By Proposition $2.1$, we can use the isomorphism from $R/I$ to $R/J$ to induce an isomorphism from $R_{m}/I_{m}$ to $R_{m}/J_{m}$ by the following: if $\phi(x_{i})=f_{i}(x_{1},...,x_{n})$, and $\nu_{m}(\phi(x_{i}))=f_{i}^{(0)}+f_{i}^{(1)}t+...f_{i}^{(m)}t^{m}$, then $x_{i}^{(j)}$ maps to $f_{i}^{(j)}$. Since $\phi$ is an isomorphism, we have $f_{i}(0)=0$ for all $i$. Noticed that $f_{i}^{(0)}(x_{1}^{(0)},...,x_{n}^{(0)})=f_{i}(x_{1}^{(0)},...,x_{n}^{(0)})$, $\mathbf{m}^{e}$ in $R/I$ will map to $\mathbf{m}^{e}$ in $R/J$. So there exists a map $\phi_{m}: R_{m}[t]/(t^{m+1},I_{m},\mathbf{m}^{e}) \longrightarrow R_{m}[t]/(t^{m+1},J_{m},\mathbf{m}^{e})$ 
such that the following diagram is commutative. 

\begin{center}
	\begin{tikzpicture}
		\draw [->](0,0) node[left]{$ker(\mu_{m})$}--(2,0)node[right]{$R$};
		\draw [->](2.5,0) --node[pos=0.5,above]{$\mu_{m}$}(4.5,0)node[right]{$R_{m}[t]/(t^{m+1},I_{m},\mathbf{m}^{e})$};
		\draw [->,dashed](-1,-0.2) --(-1,-2)node[below]{$ker(\nu_{m})$};
		\draw [->](0,-2.3) --(2,-2.3)node[right]{$R$};
		\draw [->](2.5,-2.3) --node[pos=0.5,above]{$\nu_{m}$}(4.5,-2.3)node[right]{$R_{m}[t]/(t^{m+1},J_{m},\mathbf{m}^{e})$};
		\draw [->](2.25,-0.2) --node[pos=0.5,right]{$\phi$}(2.25,-2);
		\draw [->](6,-0.2) --node[pos=0.5,right]{$\phi_{m}$}(6,-2);
	\end{tikzpicture}.
\end{center}

Because the diagram is commutative, $ker(\mu_{m})$ will be mapped to $ker(\nu_{m})$ through $\phi$. Similarly, this procedure can be done when $\phi$ is substituted by $\phi^{-1}$, so $\phi$ gives isomorphism between $I^{m-jc}$ and $J^{m-jc}$.
\end{proof}

According to Fact $2.4$, generally we have $I+\mathbf{m}^{m+1} \subset I^{m-jc}$ for any ideal $I \subset R$. Through calculating, we know some of the ideals have the property that the relation above is equality, we call them "good" jet closure:
\begin{definition}
For any ideal $I \subset R$, we say the $m$-th jet closure(resp.jet support closure) is "good" if $I+\mathbf{m}^{m+1}=I^{m-jc}(resp.I^{m-jsc})$.
\end{definition}

Next we will prove Theorem \ref{mtha}: when the ideal $I$ is extended from a homogeneous ideal in $k[x_{1},...,x_{n}]$, its every jet closure is good:
\begin{theorem}
Let $R=k[[x_{1},...,x_{n}]]$ and $I \subset R$ be an ideal. If $I$ is an ideal extended from a homogeneous ideal in $k[x_{1},...,x_{n}]$, i.e. $I=(f_{1},...,f_{s})$, where $f_{1},...,f_{s}$ are homogeneous polynomials from $k[x_{1},...,x_{n}]$, then the $m$-th jet closure of $I$ is $I^{m-jc}=I+\mathbf{m}^{m+1}$.
\end{theorem}
\begin{proof}
Assume $I=(f_1,f_2,...,f_{s})$ and $f_1,...,f_{s}$ are homogeneous polynomials with degree $d_1,...,d_{s}$. Since $I$ and $\mathbf{m}^{m+1}$ are contained in $I^{m-jc}$ according to Fact $2.15$, we have a ring map $\hat{\mu_{m}}:R/(I+\mathbf{m}^{m+1}) \longrightarrow R_{m}[t]/(t^{m+1},I_{m},\mathbf{m}^{e})$ 
induced by $\mu_{m}$. It suffices to show that $\hat{\mu_{m}}$ is injective. 

Now define a map $\alpha :R/(I+\mathbf{m}^{m+1}) \longrightarrow (R/(I+\mathbf{m}^{m+1}))[t]/(t^{m+1})$ such that $\alpha(x_{i})=x_{i}t$, then one can check that $\alpha$ is injective. In order to define a map $\beta$ from $R_{m}[t]/(t^{m+1},I_{m},\mathbf{m}^{e})$ to $(R/(I+\mathbf{m}^{m+1}))[t]/(t^{m+1})$, we let $\beta(x_{i}^{(j)})=0$ when $j \neq 1$ and $\beta(x_{i}^{(j)})=x_{i}$ when $j=1$. Also we define $\beta(t)=t$ and we need to check whether $\beta(I_{m})=0$ in $(R/(I+\mathbf{m}^{m+1}))[t]/(t^{m+1})$. Assume $f_{i}(x_{1}^{(1)}t+...+x_{1}^{(m)}t^{m},...,x_{n}^{(1)}t+...+x_{n}^{(m)}t^{m})=f_{i}^{(1)}t+...+f_{i}^{(m)}t^{m}$ in $R_{m}[t]/(t^{m+1},I_{m},\mathbf{m}^{e})$, then $I_{m}$ is defined by all the $f_{i}^{(j)}$. Because $f_{i}$ is homogeneous with degree $d_{i}$, $f_{i}^{(1)}=f_{i}^{(2)}=...=f_{i}^{(d_{i}-1)}=0$ and when $j>d_{i}$, each term in $f_{i}^{(j)}$ contains at least one $x_{p}^{(q)}$ where $q \neq 1$. It follows that $\beta(f_{i}^{(j)})=0$ when $j \neq d_{i}$, and $\beta(f_{i}^{(d_{i})})=f_{i}(x_{1},...,x_{n})$, so $\beta(I_{m})=0$ in $(R/(I+\mathbf{m}^{m+1}))[t]/(t^{m+1})$. 

\begin{center}
\
\xymatrix@=8ex
{
 R/(I+\mathbf{m}^{m+1}) \ar[r]^{\hat{\mu_{m}}} \ar[dr]^{\alpha} & R_{m}[t]/(t^{m+1},I_{m},\mathbf{m}^{e}) \ar[d]^{\beta} \\
 & R/(I+\mathbf{m}^{m+1})[t]/(t^{m+1}).
}
\end{center}

Noticed that $\alpha(x_{i})=\beta(\hat{\mu_{m}}(x_{i}))$ for all $i$, so $\alpha=\beta \circ \hat{\mu_{m}}$. Since $\alpha$ is injective, $\hat{\mu_{m}}$ is injective and the proposition follows.
\end{proof}

\begin{remark}
The jet support closure of homogeneous ideal may not be good, a simple example is $I=(x^{2},y^{2})$ in $R=k[[x,y]]$. Direct computation shows that $I^{2-jsc}=I+\mathbf{m}^{3}+(xy)$.
\end{remark}

In particular, if the ideal $I$ is defined by a homogeneous polynomial $f$, we can use the proposition above to calculate the $k$-dimension of $R/I^{m-jc}$.
\begin{corollary}
If $I=(f)$ is defined by a homogeneous polynomial $f$ with degree $d$, then $dim_{k}R/I^{m-jc}=\binom{n+m}{n}-\binom{n+m-d}{n}.$ 
\end{corollary}

\begin{proof}
From Theorem $3.1$ we know $I^{m-jc}=(f)+\mathbf{m}^{m+1}$.If $d \ge m+1$, $I^{m-jc}=\mathbf{m}^{m+1}$ and $dim_{k}R/I^{m-jc}=\binom{m+n}{n}$. If $d \le m$, $dim_{k}R/I^{m-jc}$ equals the coefficient of the $t^{m}$ of the Hilbert series $H_{(f)}(t)$. Since $H_{(f)}(t)=\frac{1-t^{d}}{(1-t)^{n+1}}$, $dim_{k}R/I^{m-jc}=\binom{n+m}{n}-\binom{n+m-d}{n}$. If we set $\binom{n+m-d}{n} < 0$ when $m < d$, then the corollary follows.
\end{proof}

Using Theorem \ref{mtha}, we can calculate the $m$-th jet closure of $A_{n}$ and $D_{n}$ singularity by hand when $m$ is small. First we compute the jet closure of $A_{n}$ singularity. Since the homogeneous case has been solved, we only need to consider about the case when $n \ge 3$. 

\begin{proposition}
Consider the $A_{n-1}$ singularity in $R=k[[x,y]]$: $I=(x^{2}+y^{n})$ $(n \ge 3)$. When $m \le n$, the $m$-th jet closure of $I$ is good, and when $m=n+1$, $I^{m-jc}=\mathbf{m}^{m+1}+I+(x^{3})$.  
\end{proposition}
\begin{proof}
The same as the proof of the proposition above, we have $I^{m-jc}=(I+\mathbf{m}^{m+1})^{m-jc}$. When $m \le n-1$, $y^{n} \in \mathbf{m}^{m+1}$, so $I^{m-jc}=((x^{2})+(y^{n})+\mathbf{m}^{m+1})^{m-jc}$. It is a monomial ideal, so by the theorem above, $I^{m-jc}=(x^{2})+(y^{n})+\mathbf{m}^{m+1}=I+\mathbf{m}^{m+1}$. 

When $m=n$, we need to show if $f \in R$ satisfies $\mu_{m}(f)=0$ in $R_{m}[t]/(t^{m+1},I_{m},\mathbf{m}^{e})$, then $f \in \mathbf{m}^{m+1}+I$. First we calculate $(t^{m+1},I_{m},\mathbf{m}^{e})$. Let $x=x_{1}t+...+x_{n}t^{n}, y=y_{1}t+...+y_{n}t^{n}$ in $x^{2}+y^{n}$, we get
\begin{flalign*}
(t^{m+1},I_{m},\mathbf{m}^{e})=(t^{m+1},x_{0},y_{0},f_{2},...,f_{n-1},f_{n}+y_{1}^{n}),	
\end{flalign*} 
where $f_{i}=\sum_{j=1}^{i-1}x_{j}x_{i-j}$. Assume $f=a_{1}y+a_{2}y^{2}+...+a_{n}y^{n}+b_{0}x+b_{1}xy+...+b_{n-1}xy^{n-1}$ and $\mu_{n}(f) \in (t^{m+1},I_{m},\mathbf{m}^{e})$. Consider the quotient map:
\begin{flalign*}
q_{1}: R_{m}[t]/(t^{m+1},I_{m},\mathbf{m}^{e}) \longrightarrow R_{m}[t]/(t^{m+1},I_{m},\mathbf{m}^{e},x_{1},...,x_{n}), 
\end{flalign*}
so $q_{1} \circ \mu_{n}(f)=0$. $R_{m}[t]/(t^{m+1},I_{m},\mathbf{m}^{e},x_{1},...,x_{n})=R_{m}[t]/(t^{m+1}, x_{0},...,x_{n},y_{0},y_{1}^{n})$, so $a_{1}=a_{2}=...=a_{n-1}=0$. Similarly consider another quotient map:
\begin{flalign*}
q_{2}:R_{m}[t]/(t^{m+1},I_{m},\mathbf{m}^{e}) \longrightarrow R_{m}[t]/(t^{m+1},I_{m},\mathbf{m}^{e},x_{2},...,x_{n}), 
\end{flalign*}
and we can get $q_{2} \circ \mu_{n}(f) \in (t^{m+1},I_{m},\mathbf{m}^{e},x_{2},...,x_{n})=(t^{m+1},x_{0},x_{2},...,x_{n},y_{1}^{n},x_{1}^{2})$, so $b_{1}=...=b_{n-1}=0$ and $f=a_{n}y^{n}$. Hence $y_{1}^{n} \in (t^{m+1},x_{0},y_{0},f_{2},...,f_{n-1},f_{n}+y_{1}^{n})$. Let $x_{0}=x_{2}=x_{3}=...=x_{n-2}=x_{n}=y_{2}=...=y_{n}=0$, we have $y_{1}^{n} \in (x_{1}^{2},y_{1}^{n}+2x_{1}x_{n-1})$, this is a contradiction.  

When $m=n+1$, we have to prove two things. First we prove that $x^{3} \in I^{(n+1)-jc}$. It is equivalent to prove $xy^{n} \in I^{(n+1)-jc}$. Through calculating, we have the following equation which uses the condition $n \ge 3$:
\begin{flalign*}
(t^{n+2},I_{n+1},\mathbf{m}^{e})=(t^{m+1},x_{0},y_{0},f_{2},...,f_{n-1},f_{n}+y_{1}^{n},f_{n+1}+ny_{1}^{n-1}y_{2}),
\end{flalign*}
where $f_{i}=\sum_{j=1}^{i-1}x_{j}x_{i-j}$. Noticed that $\mu_{n+1}(xy^{n})=x_{1}y_{1}^{n}$, it suffices to show $x_{1}f_{n} \in (f_{2},...,f_{n-1})$. Let $c_{2},...,c_{n-1} \in k$ and we want to find suitable coefficients such that $x_{1}f_{n}=c_{n-1}x_{n-1}f_{2}+...+c_{2}x_{2}f_{n-1}$. We compare the coefficient of every monomial of each side and solve the linear equation for each $c_{i}$. 

Case $1$: $n=2s+1$. 

For each monomial $x_{i}x_{j}x_{l}$ satisfying $1 \le i \le j \le l, i+j+l=2s+2$, the coefficient of it at left side is: 
\begin{flalign*}
2, when \ i=1; \\
0, when \ i>1.
\end{flalign*}
The coefficient of it at right side is: 
\begin{flalign*}
c_{2s}, when \ i=j=1; \\
2c_{j}+2c_{l}, when \ i=1,j>1; \\
c*(c_{i}+c_{j}+c_{l}), when \ i>1.
\end{flalign*}

Where $c \in \{1,2,1/3\}$ which depends on if $i,j,l$ is equal to each other. So we get the equations: 
\begin{flalign*}
c_{2s}&=2; \\
c_{j}+c_{2s+1-j}&=1;(j=2,...,2s-1); \\
c_{i}+c_{j}+c_{l}&=0, 2 \le i \le j \le l,i+j+l=2s+2.
\end{flalign*}

Let $c_{2s}=2$ and $c_{i}=(3i-2s-2)/(2s-1)(i=2,...,2s-1)$, one can check that this is the solution of the equations, so we are done. 

Case $2$: $n=2s$. 

For each monomial $x_{i}x_{j}x_{l}$ satisfying $1 \le i \le j \le l, i+j+l=2s+1$, the coefficient of it at left side is: 
\begin{flalign*}
2, when \ i=1,j<s; \\
1, when \ i=1,j=l=s; \\
0, when \ i>1.
\end{flalign*}
The coefficient of it at right side is: 
\begin{flalign*}
c_{2s-1}, when \ i=j=1; \\
2c_{j}+2c_{l}, when \ i=1,s>j>1; \\
2c_{s}, when \ i=1,j=s; \\
c*(c_{i}+c_{j}+c_{l}), when \ i>1,
\end{flalign*}
where $c \in \{1,2,1/3\}$ which depends on if $i,j,l$ is equal to each other. So we get the equations: 
\begin{flalign*}
c_{2s-1}&=2; \\
c_{j}+c_{2s-j}&=1;(j=2,...,s-1) \\
c_{s}&=1/2; \\
c_{i}+c_{j}+c_{l}&=0, 2 \le i \le j \le l,i+j+l=2s+1.
\end{flalign*}

Let $c_{2s-1}=2$ and $c_{i}=(3i-2s-1)/(2s-2)(i=2,...,2s-2)$, one can check that this is the solution of the equations, so we are done. 

The second thing we need to do, is to show $I^{(n+1)-jc} \subset I+\mathbf{m}^{n+2}+(x^{3})$. Suppose there is $f \in R$, $f=a_{1}y+...+a_{n+1}y^{n+1}+b_{0}x+b_{1}xy+...+b_{n-1}xy^{n-1}$ and $\mu_{n+1}(f)=0$. Similarly we consider the quotient map:
\begin{flalign*}
q_{1}:R_{n+1}[t]/(t^{n+2},I_{n+1},\mathbf{m}^{e}) \longrightarrow R_{n+1}[t]/(t^{n+2},I_{n+1},\mathbf{m}^{e},x_{1},...,x_{n+1}), 
\end{flalign*}
so $q_{1} \circ \mu_{n+1}(f) \in (t^{n+2},I_{n+1},\mathbf{m}^{e},x_{1},...,x_{n+1})=(t^{n+2},y_{0},y_{1}^{n},ny_{2}y_{1}^{n-1},x_{0},...,x_{n+1})$. This tells us $a_{1}=...=a_{n-1}=0$. Consider another quotient map:
\begin{flalign*}
q_{2}:R_{n+1}[t]/(t^{n+2},I_{n+1},\mathbf{m}^{e}) \longrightarrow R_{n+1}[t]/(t^{n+2},I_{n+1},\mathbf{m}^{e},x_{2},...,x_{n+1}), 
\end{flalign*}
and we can get $b_{0}=...=b_{n-1}=0$. So $f=a_{n}y^{n}+a_{n+1}y^{n+1}$ and $y_{1}^{n+1} \in (t^{n+2},I_{n+1},\mathbf{m}^{e})$. Let $x_{2}=...=x_{n-2}=x_{n}=0$ and we have $y_{1}^{n+1} \in (x_{1}^{2},y_{1}^{n}+2x_{1}x_{n-1},ny_{2}y_{1}^{n-1})$, this is a contradiction.
\end{proof}

\begin{corollary}
When $I=(x^{2}+y^{n})(n \ge 3)$ in $R=\mathbb{C}\{x,y\}$, we can calculate the dimension as follows: $dim_{k}R/I^{m-jc}=\begin{cases} 2m+1 &\quad if \ m \le n \\ 2m &\quad if \ m=n+1 \end{cases}$.
\end{corollary}
\begin{proof}
When $m \le n$, $I^{m-jc}=I+\mathbf{m}^{m+1}$, and a $k$-basis of $R/I^{m-jc}$ is $1,y,...,y^{m},x,xy,...xy^{m-1}$, so dimension is $2m+1$. When $m=n+1$, $I^{m-jc}=I+\mathbf{m}^{n+2}+(x^{3})$, and a $k$-basis of $R/I^{m-jc}$ is $1,y,...,y^{n+1},x,xy,..,xy^{n-1}$, so dimension is $2m$.
\end{proof}

\begin{proposition}
Consider the $D_{n}$ singularity: $I=(x^{2}y+y^{n-1})$. When $m \le n-1$, the $m$-th jet closure of $I$ is good.  
\end{proposition}
\begin{proof}
When $m<n-1$, the same as the $A_{n}$ case, $I^{m-jc}=(I+\mathbf{m}^{m+1})^{m-jc}=(x^{2}y)^{m-jc}$ since $y^{n-1} \in \mathbf{m}^{m+1}$. By the monomial case, $I^{m-jc}=(x^{2}y)+\mathbf{m}^{m+1}=I+\mathbf{m}^{m+1}$.

When $m=n-1$, suppose $f=a_{1}y+...+a_{n-1}y^{n-1}+b_{0}x+...+b_{n-2}xy^{n-2}+cx^{2} \in I^{(n-1)-jc}$, so $f \in ker(\mu_{n-1})$. Similarly we consider the map:
\begin{flalign*}
q_{1}:R_{n-1}[t]/(t^{n},I_{n-1},\mathbf{m}^{e}) \longrightarrow R_{n-1}[t]/(t^{n},I_{n-1},\mathbf{m}^{e},x_{2},...,x_{n-1},y_{2},...,y_{n-1}).
\end{flalign*}
We have $a_{1}=...=a_{n-2}=b_{0}=...=b_{n-2}=c=0$, i.e. $f=a_{n-1}y^{n-1}$, since $(t^{n},I_{n-1},\mathbf{m}^{e},x_{2},...,x_{n-1},y_{2},...,y_{n-1})=(t^{n},x_{0},x_{2},...,x_{n-1},y_{0},y_{2},...,y_{n-1},y_{1}x_{1}^{2},y_{1}^{n-1})$. Consider another map:
\begin{flalign*}
q_{2}:R_{n-1}[t]/(t^{n},I_{n-1},\mathbf{m}^{e}) \longrightarrow R_{n-1}[t]/(t^{n},I_{n-1},\mathbf{m}^{e},x_{2},...,x_{n-4},x_{n-2},y_{2},...,y_{n-1}) .
\end{flalign*}
The ideal at the right side is $(t^{n}, x_{0},x_{2},...,x_{n-4},x_{n-2},y_{0},y_{2},...,y_{n-1},x_{1}^{2}y_{1},y_{1}^{n-1}+2y_{1}x_{1}x_{n-3})$.
If $a_{n-1} \neq 0$, $y_{1}^{n-1} \in (x_{1}^{2}y_{1},y_{1}^{n-1}+y_{1}x_{1}x_{n-3})$, this is a contradiction.  
\end{proof}

Using the code we have given in appendix, we can calculate $m$-th jet  closure of $E_{6},E_{7},E_{8}$ singularities when $m$ is small. Combining with Proposition $3.2$ and Proposition $3.3$,  we obtain the proof of Theorem B.
\begin{theorem}
Let $R=k[[x,y]]$, $I \subset R$ defines a simple curve singularity.
	
(1)$A_{n-1}$ singularity, $I=(x^{2}+y^{n})$ $(n \ge 3)$. When $m \le n$, $I^{m-jc}=I+\mathbf{m}^{m+1}$, and when $m=n+1$, $I^{m-jc}=\mathbf{m}^{m+1}+I+(x^{3})$.  
	
(2)$D_{n}$ singularity, $I=(x^{2}y+y^{n-1})$. When $m \le n-1$, $I^{m-jc}=I+\mathbf{m}^{m+1}$.  
	
(3)$E_{6}$ singularity, $I=(x^{3}+y^{4})$. When $m \le 4$, $I^{m-jc}=I+\mathbf{m}^{m+1}$.
	
(4)$E_{7}$ singularity, $I=(x^{3}+xy^{3})$. When $m \le 4$, $I^{m-jc}=I+\mathbf{m}^{m+1}$.
	
(5)$E_{8}$ singularity, $I=(x^{3}+y^{5})$. When $m \le 5$, $I^{m-jc}=I+\mathbf{m}^{m+1}$. 
\end{theorem}

\section{The Calculation of Jet Support Closure for Monomial Ideals}
\label{s4}

In this section, we mainly focus on the calculation of jet support closure for monomial ideals. Firstly we can also prove that the jet support closure is an invariant of singularity, and the proof is similar to that of Proposition $3.1$.
\begin{proposition}
Fix $R=k[[x_1,...x_{n}]]$, $I,J \subset R$ are ideals, assume there exists an isomorphism $\phi:R \longrightarrow R$ such that $\phi(I)=J$, then $\phi$ induces $R/I^{m-jsc} \cong R/J^{m-jsc}$ for any integer $m$, $m \ge 1$. 
\end{proposition}
\begin{proof}
Similarly we assume $R_{m}=k[[x_{i}^{(j)}]](1 \le i \le n,0 \le j \le m)$, and we define the two ring maps: 
\begin{flalign*}
\lambda_{m}:R \longrightarrow R_{m}[t]/(t^{m+1},\sqrt{I_{m}+\mathbf{m}^{e}}),
\end{flalign*}
\begin{flalign*}
\Lambda_{m}:R \longrightarrow R_{m}[t]/(t^{m+1},\sqrt{J_{m}+\mathbf{m}^{e}}) .
\end{flalign*}
From $\phi: R/I \longrightarrow R/J$ we can get the map from $R_{m}/I_{m}$ to $R_{m}/J_{m}$, and from $R_{m}/\sqrt{I_{m}+\mathbf{m}^{e}}$ to $R_{m}/\sqrt{J_{m}+\mathbf{m}^{e}}$. Similarly this map gives the map:
\begin{flalign*}
\phi_{m}:R_{m}[t]/(t^{m+1},\sqrt{I_{m}+\mathbf{m}^{e}}) \longrightarrow R_{m}[t]/(t^{m+1},\sqrt{J_{m}+\mathbf{m}^{e}}) 
\end{flalign*}
and the following diagram is commutative:   
\begin{center}
	\begin{tikzpicture}
		\draw [->](0,0) node[left]{$ker(\lambda_{m})$}--(2,0)node[right]{$R$};
		\draw [->](2.5,0) --node[pos=0.5,above]{$\lambda_{m}$}(4.5,0)node[right]{$R_{m}[t]/(t^{m+1},\sqrt{I_{m}+\mathbf{m}^{e}})$};
		\draw [->,dashed](-1,-0.2) --(-1,-2)node[below]{$ker(\Lambda_{m})$};
		\draw [->](0,-2.3) --(2,-2.3)node[right]{$R$};
		\draw [->](2.5,-2.3) --node[pos=0.5,above]{$\Lambda_{m}$}(4.5,-2.3)node[right]{$R_{m}[t]/(t^{m+1},\sqrt{J_{m}+\mathbf{m}^{e}})$};
		\draw [->](2.25,-0.2) --node[pos=0.5,right]{$\phi$}(2.25,-2);
		\draw [->](6.5,-0.2) --node[pos=0.5,right]{$\phi_{m}$}(6.5,-2);
	\end{tikzpicture}.
\end{center}

Because the diagram is commutative, $ker(\lambda_{m})$ will be mapped to $ker(\Lambda_{m})$ through $\phi$. Similarly, this procedure can be done when $\phi$ is substituted by $\phi^{-1}$, so $\phi$ gives isomorphism between $I^{m-jsc}$ and $J^{m-jsc}$.
\end{proof}

If the ideal $I$ is a monomial ideal with square-free generators, we can show that the jet support closure of $I$ is good. First we give two lemmas in \cite{2} which compute $\sqrt{I_{m}}$.

\begin{lemma}
(\cite[Theorem 2.1]{2})Suppose $I \subset R$ is a monomial ideal with square-free generators, then $\sqrt{I_{m}}$ is also a square-free monomial ideal with the following generators: for each monomial generators of $I$, for example, $x_{1}x_{2}...x_{r}$, then the generators of $\sqrt{I_{m}}$ will contain the set $A_{1,2,...,r}= \{ x_{1}^{(i_1)}...x_{r}^{(i_{r})}: i_{1}+...+i_{r} \le m \}$. If the set of minimal monomial generators of $I$ is $B$, then the minimal monomial generators of $\sqrt{I_{m}}$ are $\bigcup_{x_{j_{1}}...x_{j_{r}} \in B}A_{j_{1},...,j_{r}}$.
\end{lemma}

\begin{lemma}
(\cite[Theorem 3.1]{2})Suppose $I \subset R$ is a monomial ideal, then $\sqrt{I_{m}}$ is a square-free monomial ideal with the following generators: for each monomial generators of $I$, for example, $x_{1}^{a_{1}}x_{2}^{a_{2}}...x_{r}^{a_{r}}$, then the generators of $\sqrt{I_{m}}$ will contain all the square-free monomial of the form:
\begin{flalign*}
	\sqrt{x_{1}^{(i_1)}x_{1}^{(i_{2})}...x_{1}^{(i_{a_{1}})}x_{2}^{(i_{a_{1}+1})}...x_{r}^{(i_{a_{1}+...+a_{r}})}}  \ where \ i_{1}+...+i_{a_{1}+...+a_{r}} \le m. 
\end{flalign*}
The collection of such monomials are the minimal generators of $\sqrt{I_{k}}$ if we range all the generators of I.
\end{lemma}

There is a corollary from lemma $4.2$ which will be used in the proof later.

\begin{corollary}
Assume $I \subset R$ is a monomial ideal, then in $R_{m}$, we have $\sqrt{I_{m}}+\mathbf{m}^{e}=\sqrt{I_{m}+\mathbf{m}^{e}}$. 
\end{corollary}

\begin{proof}
Since $\sqrt{I_{m}} \subset \sqrt{I_{m}+\mathbf{m}^{e}}$ and $\mathbf{m}^{e} \subset \sqrt{I_{m}+\mathbf{m}^{e}}$, we have $\sqrt{I_{m}}+\mathbf{m}^{e} \subset \sqrt{I_{m}+\mathbf{m}^{e}}$. On the other hand, $I_{m}+\mathbf{m}^{e} \subset \sqrt{I_{m}}+\mathbf{m}^{e}$, so $\sqrt{I_{m}+\mathbf{m}^{e}} \subset \sqrt{\sqrt{I_{m}}+\mathbf{m}^{e}}$, which equals to $\sqrt{I_{m}}+\mathbf{m}^{e}$ because $\sqrt{I_{m}}$ is a square-free monomial ideal.
\end{proof}

Next we give the proposition which implies the goodness of jet support closure of monomial ideal when the generators are square-free.

\begin{proposition}
If $I \subset R$ is a monomial ideal with square-free generators, then $I^{m-jsc}=I+\mathbf{m}^{m+1}$. In particular, since $(I+\mathbf{m}^{m+1}) \subset I^{m-jc} \subset I^{m-jsc}$, $I^{m-jc}=I+\mathbf{m}^{k+1}$. 
\end{proposition}
\begin{proof}
The same as in Theorem $3.1$, it suffices to show the map:
\begin{flalign*}
\hat{\lambda_{m}}:R/(I+\mathbf{m}^{m+1}) \longrightarrow R_{m}[t]/(t^{m+1},\sqrt{I_{m}},\mathbf{m}^{e}) 
\end{flalign*}
which is induced by $\lambda_{m}:R \longrightarrow R_{m}[t]/(t^{m+1},\sqrt{I_{m}},\mathbf{m}^{e})$, is injective. The definition of $\alpha :R/(I+\mathbf{m}^{m+1}) \longrightarrow (R/(I+\mathbf{m}^{m+1}))[t]/(t^{m+1})$,  $\beta : R_{m}[t]/(t^{m+1},\sqrt{I_{m}},\mathbf{m}^{e}) \longrightarrow (R/(I+\mathbf{m}^{m+1}))[t]/(t^{m+1})$ are the same as in Theorem $3.1$, so we only need to check if $\beta(\sqrt{I_{m}})$ is $0$ in $(R/(I+\mathbf{m}^{m+1}))[t]/(t^{m+1})$. Suppose the set of minimal monomial generators of $I$ is $B$, then according to Lemma $4.1$, the minimal monomial generators of $\sqrt{I_{m}}$ are $\bigcup_{x_{j_{1}}...x_{j_{r}} \in B}A_{j_{1},...,j_{r}}$. For any $x_{j_{1}}^{(i_{1})}...x_{j_{r}}^{(i_{r})}$ in $A_{j_{1},...,j_{r}}$, if $\beta(x_{j_{1}}^{(i_{1})}...x_{j_{r}}^{(i_{r})})$ is not zero, then $i_{1}=...=i_{r}=1$,thus $\beta(x_{j_{1}}^{(i_{1})}...x_{j_{r}}^{(i_{r})})=x_{j_{1}}...x_{j_{r}}$ which is zero in $(R/(I+\mathbf{m}^{m+1}))[t]/(t^{m+1})$. Similarly in Theorem $3.1$, we have $\alpha=\beta \circ \hat{\lambda_{m}}$. Since $\alpha$ is injective, $\hat{\lambda_{m}}$ is injective and the proposition follows.
\end{proof}

The calculation of the jet support closure of a general monomial ideal is not difficult, but with a lot of work. We will work by the following several steps.

\begin{proposition}
Suppose $I \subset R$ is a monomial ideal, then $I^{m-jsc}$ is also a monomial ideal for every $m \in \mathbb{N}$. 
\end{proposition}
\begin{proof}
Consider the map $\lambda_{m}:R \longrightarrow R_{m}[t]/(t^{m+1},\sqrt{I_{m}},\mathbf{m}^{e})$. 
It suffices to show that for any $f \in ker(\lambda_{m})$, every monomial in $f$ is in $ker(\lambda_{m})$. Suppose $f=\sum_{a_{1},...,a_{n}}c_{a_{1},...,a_{n}}x_{1}^{a_{1}}x_{2}^{a_{2}}...x_{n}^{a_{n}}$, and calculate $f(x_{1}^{(1)}t+...+x_{1}^{(m)}t^{m},...,x_{n}^{(1)}t+...+x_{n}^{(m)}t^{m})$. Noticed that every monomial in $R_{m}$ is of the form $\prod_{i,j} (x_{i}^{(j)})^{b_{i,j}}$, and if some coefficient of $f(x_{1}^{(1)}t+...+x_{1}^{(m)}t^{m},...,x_{n}^{(1)}t+...+x_{n}^{(m)}t^{m})$ is $\prod_{i,j} (x_{i}^{(j)})^{b_{i,j}}$, it can only be induced from $\prod_{i}x_{i}^{\sum_{j} jb_{i,j}}$. This shows that the expansion of $f(x_{1}^{(1)}t+...+x_{1}^{(m)}t^{m},...,x_{n}^{(1)}t+...+x_{n}^{(m)}t^{m})$ has no cancellation. Since $(t^{m+1},\sqrt{I_{m}},\mathbf{m}^{e})$ is a monomial ideal according to the lemma above, every monomial of $f(x_{1}^{(1)}t+...+x_{1}^{(m)}t^{m},...,x_{n}^{(1)}t+...+x_{n}^{(m)}t^{m})$ is in $(t^{m+1},\sqrt{I_{m}},\mathbf{m}^{e})$. Noted that the expansion of every monomial of $f$ is a subset of the expansion of $f$ because the expansion of $f$ has no cancellation, so every monomial of $f$ is in $ker(\lambda_{m})$.
\end{proof}

In order to compute the jet support closure of monomial ideal $I$, we can start from principal ideal.
\begin{proposition}
For a principal monomial ideal $I=(x_{1}^{a_{1}}x_{2}^{a_{2}}...x_{r}^{a_{r}})$ satisfying $r \le n, a_{1},...,a_{r} \in \mathbb{N}$, the monomials of the $m$-th jet support closure of $I$ is composed of all the following monomials: $x_{1}^{b_{1}}x_{2}^{b_{2}}...x_{n}^{b_{n}}$, where $(b_{1},b_{2},...,b_{n})$ satisfies the following property.
\begin{align}
\forall t_{1},...,t_{r} \in \mathbb{N} \  s.t. \  t_{1}b_{1}+...+t_{n}b_{n} \le m, t_{1}a_{1}+...+t_{r}a_{r} \le m.
\end{align}
\end{proposition}
\begin{proof}
Let $B=\{(b_{1},...,b_{n}):(b_{1},...,b_{n})\, \text{satsifying}\, \ (1) \}$. From the definition we can see that if $b_{1} \le c_{1},...,b_{n} \le c_{n}$ and $(b_{1},...,b_{n}) \in B$, then $(c_{1},...,c_{n}) \in B$, so $\{ x_{1}^{b_{1}}...x_{r}^{b_{n}} : (b_{1},...,b_{n}) \in B \}$ indeed comprise the monomials of some ideal $J$. 

We first show that every $x_{1}^{b_{1}}...x_{r}^{b_{n}}$ such that $(b_{1},...,b_{n}) \in B$ is in $I^{m-jsc}$. We consider every monomial of the  expansion of $(x_{1}^{(1)}t+...+x_{1}^{(m)}t^{m})^{b_{1}}...(x_{n}^{(1)}t+...+x_{n}^{(m)}t^{m})^{b_{n}}$, which is of the form $x_{1}^{(i_{1})}...x_{1}^{(i_{b_{1}})}x_{2}^{(i_{b_{1}+1})}...x_{r}^{(i_{b_{1}+...+b_{n}})}t^{i_{1}+...+i_{b_{1}+...+b_{n}}}$, where $1 \le i_{1},...,i_{b_{1}+...+b_{n}} \le m$. Without loss of generality, we may assume $i_{1} \le ... \le i_{b_{1}}$,$i_{b_{1}+1} \le ... \le i_{b_{1}+b_{2}}$,...,$i_{b_{1}+...+b_{n-1}+1} \le ... \le i_{b_{1}+...b_{n}}$. If $i_{1}+...+i_{b_{1}+...b_{n}} > m$, $x_{1}^{(i_{1})}...x_{1}^{(i_{b_{1}})}x_{2}^{(i_{b_{1}+1})}...x_{r}^{(i_{b_{1}+...+b_{n}})}t^{i_{1}+...+i_{b_{1}+...+b_{n}}} \in (t^{m+1},\sqrt{I_{m}},\mathbf{m}^{e})$. If $i_{1}+...+i_{b_{1}+...b_{n}} \le m$, $i_{1}b_{1}+...+i_{b_{1}+...+b_{n-1}+1}b_{n} \le m$, and $i_{1}a_{1}+...+i_{b_{1}+...+b_{r-1}+1}a_{r} \le m$ because $(b_{1},...,b_{n}) \in B$. By Lemma $4.2$ we know $\sqrt{(x_{1}^{(i_{1})})^{a_{1}}...(x_{r}^{(i_{b_{1}+...+b_{r-1}+1})})^{a_{r}}}=x_{1}^{(i_{1})}...x_{r}^{(i_{b_{1}+...+b_{r-1}+1})} \in \sqrt{I_{m}}$, so $x_{1}^{(i_{1})}...x_{1}^{(i_{b_{1}})}x_{2}^{(i_{b_{1}+1})}...x_{n}^{(i_{b_{1}+...+b_{n}})}t^{i_{1}+...+i_{b_{1}+...+b_{n}}} \in (t^{m+1},\sqrt{I_{m}},\mathbf{m}^{e})$.  

On the other hand, for every monomial $x_{1}^{b_{1}}...x_{n}^{b_{n}} \in I^{m-jsc}$, we want to show that $(b_{1},...,b_{n}) \in B$. Suppose $t_{1},...,t_{n} \in \mathbb{N}$ satisfy $t_{1}b_{1}+...+t_{n}b_{n} \le m$. Consider a monomial $$(x_{1}^{(t_{1})})^{b_{1}}...(x_{n}^{(t_{n})})^{b_{n}}t^{t_{1}b_{1}+...+t_{n}b_{n}}$$ which is in the expansion of $(x_{1}^{(1)}t+...+x_{1}^{(m)}t^{m})^{b_{1}}...(x_{n}^{(1)}t+...+x_{n}^{(m)}t^{m})^{b_{n}}$. Since $x_{1}^{b_{1}}...x_{n}^{b_{n}} \in I^{m-jsc}$, $(x_{1}^{(t_{1})})^{b_{1}}...(x_{n}^{(t_{n})})^{b_{n}}t^{t_{1}b_{1}+...+t_{n}b_{n}} \in (t^{m+1},\sqrt{I_{m}},\mathbf{m}^{e})$, combining the fact that $t_{1}b_{1}+...+t_{n}b_{n} \le m$, we have $(x_{1}^{(t_{1})})^{b_{1}}...(x_{r}^{(t_{r})})^{b_{r}} \in \sqrt{I_{m}}$. By Lemma $4.2$ we can get $t_{1}a_{1}+...+t_{r}a_{r} \le m$.
\end{proof}

When it comes to the general case, the set of generators is more complicated with the following Theorem \ref{mthc}.
\begin{theorem}
For a monomial ideal $I \subset R$, the m-th jet support closure of I is composed of the following monomials: $x_{1}^{b_{1}}...x_{n}^{b_{n}} \in I^{m-jsc}$ if and only if for any $t_{1},...,t_{n} \in \mathbb{N}$ satisfying $t_{1}b_{1}+t_{2}b_{2}+...+t_{n}b_{n} \le m$, there exists $x_{i_{1}}^{a_{1}}...x_{i_{s}}^{a_{s}} \in I$ such that $t_{i_{1}}a_{1}+...+t_{i_{s}}a_{s} \le m$. 
\end{theorem}
\begin{proof}
The proof is the same as the proposition above, we only need to notice that the generators of $\sqrt{I_{m}}$ is the union of all $\sqrt{(x_{i_{1}}^{a_{1}}...x_{i_{s}}^{a_{s}})_{m}}$ for $x_{i_{1}}^{a_{1}}...x_{i_{s}}^{a_{s}} \in I$.
\end{proof}

\begin{remark}
It is not generally true that $I^{m-jsc}+J^{m-jsc}=(I+J)^{m-jsc}$, i.e. we can not just combine every jet support closure of principal ideals as the whole jet support closure. For example, take $R=k[[x,y]]$, $I=(x^{2})$ and $J=(y^{2})$, direct computation shows that $I^{3-jsc}=(x^{2},xy^{2},y^{4})$, $J^{3-jsc}=(y^{2},yx^{2},x^{4})$, while $(I+J)^{3-jsc}=(x^{2},xy,y^{2})$.
\end{remark}

\section{The Calculation of Jet Support Closure for Homogeneous and Weighted Homogeneous Polynomials}
\label{s5}

In this section, we calculate several special cases of jet support closures. First we prove that the jet support closure of a reduced homogeneous polynomial is good.
\begin{proposition}
Let $f(x_{1},...,x_{n})$ be a homogeneous polynomial of degree $d$, suppose $f$ is reduced, then $(f)^{m-jsc}=(f)+\mathbf{m}^{m+1}$ in $k[[x_{1},...,x_{n}]]$.  
\end{proposition}
\begin{proof}
We denote an ideal $I=(f)$ for convenience. When $m<d$, we have $f(x_{1}^{(1)}t+...+x_{1}^{(m)}t^{m},...,x_{n}^{(1)}t+...+x_{n}^{(m)}t^{m})=0(mod \ t^{m})$, so $\sqrt{I_{m}+\mathbf{m}^{e}}=\mathbf{m}^{e}$. This implies that $(f)^{m-jsc}=\mathbf{m}^{m+1}=(f)+\mathbf{m}^{m+1}$. When $m \ge d$, $I_{m}+\mathbf{m}^{e}=(x_{1}^{(0)},...,x_{n}^{(0)},F_{d},...,F_{m})$, where $F_{d},...,F_{m}$ are defined by $f(x_{1}^{(1)}t+...+x_{1}^{(m)}t^{m},...,x_{n}^{(1)}t+...+x_{n}^{(m)}t^{m})=F_{d}t^{d}+...+F_{m}t^{m}(mod \ t^{m+1})$. Noticed that $f(x_{1}^{(1)}t+...+x_{1}^{(m)}t^{m},...,x_{n}^{(1)}t+...+x_{n}^{(m)}t^{m})=t^{d}f(x_{1}^{(1)}+...+x_{1}^{(m)}t^{m-1},...,x_{n}^{(1)}+...+x_{n}^{(m)}t^{m-1})$, $F_{d},...,F_{m}$ are the generators for $I_{m-d}$ after a change of variables. If there is a polynomial $g \notin (f)+\mathbf{m}^{m+1}$ and $g \in I^{m-jsc}$, we may assume every monomial of $g$ has degree less than $m+1$. Suppose $g(x_{1}^{(1)}t+...+x_{1}^{(m)}t^{m},...,x_{n}^{(1)}t+...+x_{n}^{(m)}t^{m})=G_{1}t^{1}+...+G_{m}t^{m}(mod \ t^{m+1})$ and $G_{1},..,G_{m} \in \sqrt{(F_{d},...,F_{m})}$. First let $x_{i}^{(j)}=0(1 \le i \le n, 2 \le j \le m)$, we have $g(x_{1}^{(1)}t,...,x_{n}^{(1)}t)=G_{1}(x_{1}^{(1)},...,x_{n}^{(1)})t+...+G_{m}(x_{1}^{(1)},...,x_{n}^{(1)})t^{m}$. After this, we let $t=1$ and get $g(x_{1}^{(1)},...,x_{n}^{(1)})=G_{1}(x_{1}^{(1)},...,x_{n}^{(1)})+...+G_{m}(x_{1}^{(1)},...,x_{n}^{(1)}) \in \sqrt{F_{d}}$. Noticed that $F_{d}=f$ and it has no multiple factors, we have $g \in (f)$, this is a contradiction. This gives us the assertion.
\end{proof}

We can compute the jet support closure of a reduced polynomial $f$ when the number of variables is two. Here is the basic setting: Suppose $f \in R=k[[x,y]]$ is a weighted homogeneous polynomial, the weight of $x$ is $a$ and the weight of $y$ is $b$, satisfying $a,b \in \mathbb{N}$, $a \neq b$ and $gcd(a,b)=1$. We assume $f(x,y)=\sum_{l=1}^{s}c_{l}x^{i_{l}}y^{j_{l}}$, where $ai_{l}+bj_{l}=d$ and $c_{l} \in k/\{0\} $ for every $l$.  
\begin{proposition}
Under this setting, when $m<d$, let $A_{m}:=\{ (u,v) \in \mathbb{N}^{2}:  ui_{l}+vj_{l}\ge m+1, l=1,...,s \}$, then $(f)^{m-jsc}$ is generated by all $x^{p}y^{q}$, where $(p,q)$ satisfies $pu+qv \ge m+1$ for all $(u,v) \in A_{m}$. 
\end{proposition}
\begin{proof}
For convenience, we set $I=(f)$. Without loss of generality, we can suppose that $a>b$, $i_{1}<...<i_{s}$ and $j_{1}>...>j_{s}$. Since $a>b$ and $ai_{l}+bj_{l}=d$, we have $i_{1}+j_{1}>...>i_{s}+j_{s}$.  

First we will calculate $\sqrt{I_{m}+\mathbf{m}^{e}}$. We claim that $$\sqrt{I_{m}+\mathbf{m}^{e}}=\bigcap_{(u,v)\in A_{m}}(x_{\{u-1\}},y_{\{v-1\}}),$$where $(x_{\{u-1\}},y_{\{v-1\}})$ denote the ideal $(x_{0},...,x_{u-1},y_{0},...,y_{v-1})$ and we prove this claim by induction. When $m<i_{s}+j_{s}$, we have $f(x_{1}t+...+x_{m}t^{m},y_{1}t+...+y_{m}t^{m})=0(mod \ t^{m+1})$, so $\sqrt{I_{m}+\mathbf{m}^{e}}=\mathbf{m}^{e}$. Noticed that $A_{m}=\mathbb{N}^{2}$, and the assertion follows. Assume when $m<r$ the claim holds and $r<d$, we consider the case when $m=r$. To calculate $\sqrt{I_{r}+\mathbf{m}^{e}}$, suppose $f(x_{1}t+...+x_{r}t^{r},y_{1}t+...+y_{r}t^{r})=F_{1}t+...+F_{r}t^{r}(mod \ t^{r+1})$, and we have $\sqrt{I_{r}+\mathbf{m}^{e}}=(x_{0},y_{0})+\sqrt{(F_{1},...,F_{r})}$. Noticed that $x_{r}$ and $y_{r}$ do not appear in $F_{1},...,F_{r-1}$, so  $f(x_{1}t+...+x_{r-1}t^{r-1},y_{1}t+...+y_{r-1}t^{r-1})=F_{1}t+...+F_{r-1}t^{r-1}(mod \ t^{r})$.  This tells us $\sqrt{I_{r-1}+\mathbf{m}^{e}}=(x_{0},y_{0})+\sqrt{(F_{1},...,F_{r-1})}$. Now we claim that $\sqrt{I_{r}+\mathbf{m}^{e}}=\sqrt{\sqrt{I_{r-1}+\mathbf{m}^{e}}+(F_{r})}$. Since $I_{r}+\mathbf{m}^{e}=I_{r-1}+\mathbf{m}^{e}+(F_{r})$, the left side is contained in the right side. Also both $\sqrt{I_{r-1}+\mathbf{m}^{e}}$ and $F_{r}$ are contained in $\sqrt{I_{r}+\mathbf{m}^{e}}$, so the right side is contained in the left side. Now we use the induction hypothesis, and we can calculate as follows:
\begin{flalign*}
&\sqrt{I_{r}+\mathbf{m}^{e}} \\
=&\sqrt{\sqrt{I_{r-1}+\mathbf{m}^{e}}+(F_{r})}  \\
=&\sqrt{\bigcap_{(u,v) \in A_{r-1}}(x_{\{u-1\}},y_{\{v-1\}})+(F_{r})}.
\end{flalign*}

In order to calculate $(x_{\{u-1\}},y_{\{v-1\}})+(F_{r})$, we can expand $f(x_{u}t^{u}+...+x_{r}t^{r},...,y_{v}t^{v}+...+y_{r}t^{r})$ and see the $t^{r}$ term. Since $(u,v) \in A_{r-1}$, $ui_{l}+vj_{l} \ge r$ for all $l$, the coefficient of $t,...,t^{r-1}$ of this expansion are zero. Actually we have $f(x_{u}t^{u}+...+x_{r}t^{r},...,y_{v}t^{v}+...+y_{r}t^{r})=\sum_{i_{l}u+j_{l}v=r}t^{r}c_{l}x_{u}^{i_{l}}y_{v}^{j_{l}}(mod \ t^{r+1})$. Now we claim that this expansion has at most one term. If not, let $(i_{l_{1}},j_{l_{1}}),(i_{l_{2}},j_{l_{2}})$ be two different pairs satisfying $ui_{l_{1}}+vj_{l_{1}}=ui_{l_{2}}+vj_{l_{2}}=r$, then $u(i_{l_{1}}-i_{l_{2}})=v(j_{l_{1}}-j_{l_{2}})$. Noticed that $ai_{l_{1}}+bj_{l_{1}}=ai_{l_{2}}+bj_{l_{2}}=d$, so we have $ub=va$. Since $gcd(a,b)=1$, $u$ must be divided by $a$ and $v$ must be divided by $b$. Assume $u=ha$ and $v=hb(h \in \mathbb{N})$, then $r=ui_{l_{1}}+vj_{l_{1}}=ui_{l_{2}}+vj_{l_{2}}=hd$. This is a contradiction because we have supposed $r<d$, so the claim follows.  

Now we return to the calculation of $(x_{\{u-1\}},y_{\{v-1\}})+(F_{r})$. If $(u,v) \in A_{r}$, then there is no $l$ such that $i_{l}u+j_{l}v=r$, so $F_{r}=0$ and $(x_{\{u-1\}},y_{\{v-1\}})+(F_{r})=(x_{\{u-1\}},y_{\{v-1\}})$. If $(u,v) \notin A_{r}$, then there exists a $l$ such that $i_{l}u+j_{l}v=r$, and $(x_{\{u-1\}},y_{\{v-1\}})+(F_{r})=(x_{\{u-1\}},y_{\{v-1\}},x_{u}^{i_{l}}y_{v}^{j_{l}})$. Thus we have:
\begin{flalign*}
&\sqrt{I_{r}+\mathbf{m}^{e}} \\
=&\sqrt{\bigcap_{(u,v) \in A_{r-1}}(x_{\{u-1\}},y_{\{v-1\}})+(F_{r})} \\
=&\sqrt{\bigcap_{(u,v) \in A_{r-1}}\sqrt{(x_{\{u-1\}},y_{\{v-1\}})+(F_{r})}}  \\
=&\sqrt{(\bigcap_{(u,v) \in A_{r}}(x_{\{u-1\}},y_{\{v-1\}}))\cap (\bigcap_{ui_{l}+vj_{l}=r}(x_{\{u-1\}},y_{\{v-1\}})+\sqrt{(x_{u}^{i_{l}}y_{v}^{j_{l}})})}.
\end{flalign*}

We claim that $\bigcap_{(u,v) \in A_{r}}(x_{\{u-1\}},y_{\{v-1\}})$ is in $(x_{\{u-1\}},y_{\{v-1\}})+\sqrt{(x_{u}^{i_{l}}y_{v}^{j_{l}})}$ for every $(u,v)$ and $(i_{l},j_{l})$ satisfying $ui_{l}+vj_{l}=r$. Actually, if both $i_{l},j_{l} \neq 0$, then the sum equals to $(x_{0},...,x_{u},y_{0},...,y_{v-1}) \cap (x_{0},...,x_{u-1},y_{0},...,y_{v})$. Since $i_{l},j_{l} \ge 1$, we have $(u+1)i_{l}+vj_{l} \ge r+1$ and $ui_{l}+(v+1)j_{l} \ge r+1$. Noted that for $l' \neq l$, $ui_{l'}+vj_{l'} \ge r+1$, so $(u+1,v) \in A_{r}$ and $(u,v+1) \in A_{r}$. If $j_{l}=0$ while $i_{l} \neq 0$, then the sum equals to $(x_{0},...,x_{u},y_{0},...,y_{v-1})$, and similarly we have $(u+1,v) \in A_{r}$. Again if $i_{l}=0$ we can do the same thing. This implies $\bigcap_{(u,v) \in A_{r}}(x_{\{u-1\}},y_{\{v-1\}})$ is contained in $(x_{\{u-1\}},y_{\{v-1\}})+\sqrt{(x_{u}^{i_{l}}y_{v}^{j_{l}})}$. Finally we can compute $\sqrt{I_{r}+\mathbf{m}^{e}}$ as follows:
\begin{flalign*}
&\sqrt{I_{r}+\mathbf{m}^{e}} \\
=&\sqrt{(\bigcap_{(u,v) \in A_{r}}(x_{\{u-1\}},y_{\{v-1\}}))\cap (\bigcap_{ui_{l}+vj_{l}=r}(x_{\{u-1\}},y_{\{v-1\}})+\sqrt{(x_{u}^{i_{l}}y_{v}^{j_{l}})})} \\
=&\sqrt{\bigcap_{(u,v) \in A_{r}}(x_{\{u-1\}},y_{\{v-1\}})} \\
=&\bigcap_{(u,v) \in A_{r}}(x_{\{u-1\}},y_{\{v-1\}}).
\end{flalign*}

So we have finished the computation of $\sqrt{I_{m}+\mathbf{m}^{e}}$ by induction.  

By Proposition $2.5 \ (3)$, let $\lambda_{m}: R \longrightarrow R_{m}[t]/(t^{m+1},\sqrt{I_{m}+\mathbf{m}^{e}})$ be the ring map defined in that proposition, then we have $I^{m-jsc}=ker(\lambda_{m})$. Since $$\sqrt{I_{m}+\mathbf{m}^{e}}=\bigcap_{(u,v)\in A_{m}}(x_{\{u-1\}},y_{\{v-1\}}),$$ $$ I^{m-jsc}=\bigcap_{(u,v)\in A_{m}}(\lambda_{m}^{-1}(t^{m+1},x_{\{u-1\}},y_{\{v-1\}})),$$ by direct computation, we can show that $$\lambda_{m}^{-1}(t^{m+1},x_{\{u-1\}},y_{\{v-1\}})=(x^{p}y^{q})(up+yq \ge m+1),$$  therefore the proposition follows.
\end{proof}

\begin{proposition}
Under the setting as the proposition above, and assume that $f$ is reduced, then when $m \ge d$, $I^{m-jsc}=(f,x^{p_{1}}y^{q_{1}})\cap(x^{p_{2}}y^{q_{2}})$, where $ap_{1}+bq_{1} \ge m+1$ and $up_{2}+vq_{2} \ge m+1$ for all $(u,v) \in A_{m}, u \le a \ or \ v \le b$. 
\end{proposition}
\begin{proof}
According to last proposition, we have $$\sqrt{I_{d-1}+\mathbf{m}^{e}}=\bigcap_{(u,v)\in A_{r-1}}(x_{\{u-1\}},y_{\{v-1\}}).$$  Similarly we can compute $\sqrt{I_{d}+\mathbf{m}^{e}}=\sqrt{\bigcap_{(u,v) \in A_{d-1}}(x_{\{u-1\}},y_{\{v-1\}})+(F_{d})}$. When $(u,v) \in A_{d}$, we have $ui_{l}+vj_{l} \ge d+1$ for every $l$, so $F_{d}=0(mod \ t^{d+1})$ and   $(x_{\{u-1\}},y_{\{v-1\}})+(F_{d})=(x_{\{u-1\}},y_{\{v-1\}})$. But when $(u,v) \notin A_{d}$, things become different from last time. If there are two different terms in the expansion of $f(x_{u}t^{u}+...+x_{d}t^{d},y_{v}t^{v}+...+y_{d}t^{d})(mod \ t^{d+1})$, namely, there are two different pairs $(i_{l_{1}},j_{l_{1}}),(i_{l_{2}},j_{l_{2}})$ such that $ui_{l_{1}}+vj_{l_{1}}=d=ui_{l_{2}}+vj_{l_{2}}$. Since $ai_{l_{1}}+bj_{l_{1}}=d=ai_{l_{2}}+bj_{l_{2}}$, we have $ub=va$, and thus $u=a,v=b$ because $gcd(a,b)=1$. This tells us things become different only when $(u,v)=(a,b)$. In particular, when $(u,v)=(a,b)$, $f(x_{u}t^{u}+...+x_{d}t^{d},y_{v}t^{v}+...+y_{d}t^{d})=t^{d}f(x_{a},y_{b})(mod \ t^{d+1})$. This tells us:
\begin{flalign*}
&\sqrt{I_{d}+\mathbf{m}^{e}} \\
=&\sqrt{\bigcap_{(u,v) \in A_{d}}(x_{\{u-1\}},y_{\{v-1\}})\cap((x_{\{a-1\}},y_{\{b-1\}})+\sqrt{f(x_{a},y_{b})})} \\
=&\sqrt{\bigcap_{(u,v) \in A_{d}}(x_{\{u-1\}},y_{\{v-1\}})\cap(x_{\{a-1\}},y_{\{b-1\}},f(x_{a},y_{b}))} \\
=&\bigcap_{(u,v) \in A_{d}}(x_{\{u-1\}},y_{\{v-1\}})\cap(x_{\{a-1\}},y_{\{b-1\}},f(x_{a},y_{b})). 
\end{flalign*}

The last equation holds because the intersection of radical ideals is still radical ideal. Noticed that when $ u > a, v > b$, $(x_{\{a-1\}},y_{\{b-1\}},f(x_{a},y_{b}))$ is contained in $(x_{\{u-1\}},y_{\{v-1\}})$, thus we can change the range of $(u,v)$ to $(u,v) \in A_{d}, u \le a \ or \ v \le b$. Now we claim when $m \ge d$, $\sqrt{I_{m}+\mathbf{m}^{e}}=\bigcap_{(u,v) \in A_{m},u \le a \ or \ v \le b}(x_{\{u-1\}},y_{\{v-1\}})\cap((x_{\{a-1\}},y_{\{b-1\}})+\sqrt{(F_{d},...,F_{m})})$, where $F_{i}$ are defined by $f(x_{a}t^{a}+...+x_{m}t^{m},y_{b}t^{b}+...+y_{m}t^{m})=F_{d}t^{d}+...+F_{m}t^{m}$. We prove it by induction. When $m=d$, we have proved the claim. Now we assume that when $m<r$ the claim holds, and we consider the case when $m=r$. The same as the last proposition, we can do the similar calculation as follows:
\begin{flalign*}
&\sqrt{I_{r}+\mathbf{m}^{e}} \\
=&\sqrt{\sqrt{I_{r-1}+\mathbf{m}^{e}}+(F_{r})}  \\
=&\sqrt{\bigcap_{\substack{(u,v) \in A_{r-1} \\ u \le a \ or \ v \le b}}(x_{\{u-1\}},y_{\{v-1\}})+(F_{r})\cap ((x_{\{a-1\}},y_{\{b-1\}})+\sqrt{(F_{d},...,F_{r-1})}+(F_{r}))} \\
=&\sqrt{\bigcap_{\substack{(u,v) \in A_{r-1} \\ u \le a \ or \ v \le b}}\sqrt{(x_{\{u-1\}},y_{\{v-1\}},F_{r})}\cap ((x_{\{a-1\}},y_{\{b-1\}})+\sqrt{(F_{d},...,F_{r})})}.
\end{flalign*}

The same as the last proposition, if $(u,v) \in A_{r}$, $F_{r}=0$, and if $(u,v) \notin A_{r}$, the expansion of $f(x_{u}t^{u}+...+x_{r}t^{r},y_{v}t^{v}+...+y_{r}t^{r})(mod \ t^{r+1})$ has only one term except $(u,v)$ can be divided by $(a,b)$. But now $u \le a$ or $v \le b$, so if the expansion above has more than one term, $(u,v)=(a,b)$, contradict to the fact that $r>d$. So we can continue to calculate:
\begin{flalign*}
&\sqrt{I_{r}+\mathbf{m}^{e}} \\
=&\sqrt{\bigcap_{\substack{(u,v) \in A_{r-1} \\ u \le a \ or \ v \le b}}\sqrt{(x_{\{u-1\}},y_{\{v-1\}},F_{r})}\cap ((x_{\{a-1\}},y_{\{b-1\}})+\sqrt{(F_{d},...,F_{r})})} \\
=&\sqrt{\bigcap_{\substack{(u,v) \in A_{r} \\ u \le a \ or \ v \le b}}(x_{\{u-1\}},y_{\{v-1\}})\cap ((x_{\{a-1\}},y_{\{b-1\}})+\sqrt{(F_{d},...,F_{r})})} \\
=&\bigcap_{\substack{(u,v) \in A_{r} \\ u \le a \ or \ v \le b}}(x_{\{u-1\}},y_{\{v-1\}})\cap ((x_{\{a-1\}},y_{\{b-1\}})+\sqrt{(F_{d},...,F_{r})}).
\end{flalign*}

By Proposition $2.5 \ (3)$, let $\lambda_{m}: R \longrightarrow R_{m}[t]/(t^{m+1},\sqrt{I_{m}+\mathbf{m}^{e}})$ be the ring map defined in that proposition, then we have $I^{m-jsc}=ker(\lambda_{m})$. Through direct computation, we can show that  $\lambda_{m}^{-1}(t^{m+1},x_{\{u-1\}},y_{\{v-1\}})=(x^{p}y^{q})(up+vq \ge m+1)$, so we only need to compute $\lambda_{m}^{-1}((x_{\{a-1\}},y_{\{b-1\}})+\sqrt{(F_{d},...,F_{m})})$.  

Now we claim: $\lambda_{m}^{-1}((x_{\{a-1\}},y_{\{b-1\}})+\sqrt{(F_{d},...,F_{m})})=(f,x^{p}y^{q})(ap+bq \ge m+1)$. Since $f(x_{a}t^{a}+...+x_{m}t^{m},y_{b}t^{b}+...+y_{m}t^{m})=t^{d}F_{d}+...+t^{m}F_{m}$ and  $(x_{a}t^{a}+...+x_{m}t^{m})^{p}(y_{b}t^{b}+...+y_{m}t^{m})^{q}=0(mod \ t^{m+1})$, so the right side is contained in the left side. To prove the other containing relation, we assume the opposite case, namely, there is a $g \in \lambda_{m}^{-1}((x_{\{a-1\}},y_{\{b-1\}})+\sqrt{(F_{d},...,F_{m})})$ but $g \notin  (f,x^{p}y^{q})(ap+bq \ge m+1)$. For every monomial $x^{i}y^{j}$ in $g$, we can assume that $ai+bj \le m$(otherwise we can delete it from $g$). Set $g(x_{a}t^{a}+...+x_{m}t^{m},y_{b}t^{b}+...+y_{m}t^{m})=tG_{1}+...+t^{m}G_{m}(mod \ t^{m+1})$, and $G_{1},...,G_{m} \in \sqrt{(F_{d},...,F_{m})}$. Let $x_{a+1}=...=x_{m}=y_{b+1}=...=y_{m}=0$, and we can get $g(x_{a}t^{a},y_{b}t^{b})=tG_{1}(x_{a},y_{b})+...+t^{m}G_{m}(x_{a},y_{b})$. This is an equation other than an equation in the module sense because we have assumed that every monomial $x^{i}y^{j}$ in $g$ satisfies that $ia+jb \le m$. Now let $t=1$ and we can get $g(x_{a},y_{b})=G_{1}(x_{a},y_{b})+...+G_{m}(x_{a},y_{b})$. Noticed that when $x_{a+1}=...=x_{m}=y_{b+1}=...=y_{m}=0$, $\sqrt{(F_{d},...,F_{m})}=\sqrt{f(x_{a},y_{b})}$, we have $g \in (f)$ because $f$ is reduced. In the last, we have $I^{m-jsc}=ker(\lambda_{m})=(f,x^{p_{1}}y^{q_{1}})\cap(x^{p_{2}}y^{q_{2}})$, where $ap_{1}+bq_{1} \ge m+1$ and $up_{2}+vq_{2} \ge m+1$ for all $(u,v) \in A_{m}, u \le a \ or \ v \le b$.  
\end{proof}

Combining with the three propositions above, we obtain the proof of the first two parts of Theorem \ref{mthd}.
\begin{theorem}\label{th5.4}
(1)Let $f(x_{1},...,x_{n})$ be a homogeneous polynomial of degree $d$, suppose $f$ is reduced, then $(f)^{m-jsc}=(f)+\mathbf{m}^{m+1}$ in $k[[x_{1},...,x_{n}]]$. 

(2)Suppose $f \in R=k[[x,y]]$ is a weighted homogeneous polynomial, the weight of $x$ is $a$ and the weight of $y$ is $b$, satisfying $a,b \in \mathbb{N}$ and $gcd(a,b)=1$. We assume $f(x,y)=\sum_{l=1}^{s}c_{l}x^{i_{l}}y^{j_{l}}$, where $ai_{l}+bj_{l}=d$ and $c_{l} \in k/\{0\} $ for every $l$. Let $A_{m}:=\{ (u,v) \in \mathbb{N}^{2}: ui_{l}+vj_{l}\ge m+1, l=1,...,s \}$, then when $m<d$, $(f)^{m-jsc}$ is generated by all $x^{p}y^{q}$, where $(p,q)$ satisfies $pu+qv \ge m+1$ for all $(u,v) \in A_{m}$, and when $m \ge d$, $I^{m-jsc}=(f,x^{p_{1}}y^{q_{1}})\cap(x^{p_{2}}y^{q_{2}})$, where $ap_{1}+bq_{1} \ge m+1$ and $up_{2}+vq_{2} \ge m+1$ for all $(u,v) \in A_{m}, u \le a \ or \ v \le b$. 
\end{theorem}

Using this Theorem \ref{th5.4}, we can calculate the jet support closure of simple singularities when the number of variables is two. Also in concrete cases we can calculate the $k$-dimension of $R/(f)^{m-jsc}$. Here we list the following results:
\begin{corollary}
In the following cases we assume ideal $I=(f)$, where $f \in R=k[[x,y]]$ defines a simple curve singularity. 

(1)Let $f=x^{2}+y^{2n}$ be the $A_{2n-1}$ singularity, then when $m<2n$,$I^{m-jsc}=(x^{2},xy^{\lceil \frac{m}{2} \rceil},y^{m+1})$, $dim_{k}R/(f)^{m-jsc}=m+1+\lceil \frac{m}{2} \rceil$; and when $m \ge 2n$, $I^{m-jsc}=(f,x^{i}y^{j})(ni+j \ge m+1)$, $dim_{k}R/(f)^{m-jsc}=2m+2-n$.    

(2)Let $f=x^{2}+y^{2n+1}$ be the $A_{2n}$ singularity, then when $m<2n+1$,$I^{m-jsc}=(x^{2},xy^{\lceil \frac{m}{2} \rceil},y^{m+1})$, $dim_{k}R/(f)^{m-jsc}=m+1+\lceil \frac{m+1}{2} \rceil$; when $2n+1 \le m <4n+2$, $I^{m-jsc}=(x^{2},xy^{\lceil \frac{m}{4} \rceil},y^{\lceil \frac{m+1}{2} \rceil})$, $dim_{k}R/(f)^{m-jsc}=\lceil \frac{3(m+1)}{4} \rceil$; and when $m \ge 4n+2$, $I^{m-jsc}=(f,x^{i}y^{j})((2n+1)i+2j \ge m+1)$, $dim_{k}R/(f)^{m-jsc}=m+1-n$.   

(3)Let $f=x^{2}y+y^{2n-1}$ be the $D_{2n}$ singularity, suppose $q,r \in \mathbb{N}$ satisfying $m=(n-1)q+r, 0 \le r \le n-2$, then when $m<2n-1$,$I^{m-jsc}=(x^{m+1},x^{2}y,xy^{\lceil \frac{m+1}{2} \rceil},y^{m+1})$, $dim_{k}R/(f)^{m-jsc}=2m+\lceil \frac{m}{2} \rceil$; and when $m \ge 2n-1$, we have $I^{m-jsc}=(x^{2}y+y^{2n-1},y^{m+1},xy^{m+2-n},...,x^{q}y^{r+1},x^{m+1})$, $dim_{k}R/(f)^{m-jsc}=3m+2-n$.  

(4)Let $f=x^{2}y+y^{2n}$ be the $D_{2n+1}$ singularity, then we have:
when $m<2n$, $I^{m-jsc}=(x^{m+1},x^{2}y,xy^{\lceil \frac{m+1}{2} \rceil},y^{m+1})$, $dim_{k}R/(f)^{m-jsc}=2m+\lceil \frac{m+1}{2} \rceil$; when $2n \le m \le 4n-1$, $I^{m-jsc}=(x^{m+1},x^{2}y,xy^{\lceil \frac{m+2}{4} \rceil},y^{\lceil \frac{m+1}{2} \rceil})$, $dim_{k}R/(f)^{m-jsc}=\lceil \frac{7m+1}{4} \rceil$ and when $m \ge 4n$,  $I^{m-jsc}=(x^{m+1},x^{2}y,xy^{n+1},y^{m+1})\cap (f,x^{i}y^{j})((2n-1)i+2j \ge m+1)$, $dim_{k}R/(f)^{m-jsc}=2m-n+1$.  

(5)Let $f=x^{3}+y^{4}$ be the $E_{6}$ singularity, then when $m \ge 12$,  $I^{m-jsc}= (f,x^{i}y^{j})((4i+3j \ge m+1)$, $dim_{k}R/(f)^{m-jsc}=m-2$.  

(6)Let $f=x^{3}+xy^{3}$ be the $E_{7}$ singularity, then when $m \ge 9$,  $I^{m-jsc}= (f,x^{i}y^{j})((3i+2j \ge m+1)$, $dim_{k}R/(f)^{m-jsc}=\lceil \frac{3m-5}{2} \rceil$.  

(7)Let $f=x^{3}+y^{5}$ be the $E_{8}$ singularity, then when $m \ge 15$,  $I^{m-jsc}= (f,x^{i}y^{j})(5i+3j \ge m+1)$, $dim_{k}R/(f)^{m-jsc}=m-3$. 
\end{corollary}

For any two simple curve singularities $R/I_{1},R/I_{2}$, let $M:=max\{\mu(R/I_{1}),\mu(R/I_{1})\}$ be the maximum  of Milnor numbers, then we can distinguish the simple curve singularities by just looking at their $k$-dimensions of their jet support closures of order not bigger than $M+1$. The concrete statement is as the third part of Theorem \ref{mthd}:

\begin{theorem}
Suppose $R/I_{1}$ and $R/I_{2}$ are two simple curve singularities. Set $M=max\{\mu(R/I_{1}),\mu(R/I_{2})\}$. Then we have $R/I_{1} \cong R/I_{2}$ if and only if $R/I_{1}^{m-jsc} \cong R/I_{2}^{m-jsc}$ for all $m \le M+1$.
\end{theorem}
\begin{proof}
We can get this conclusion by comparing the $k$-dimension of $R/I_{i}^{m-jsc}(i=1,2,m \le M+1)$. Corollary $5.1$ has given the data of case when $R/I_{i}$ is $A_{n},D_{n}$ singularity. We will give the jet support closures for small orders of $E_{6},E_{7},E_{8}$ curve singularities in the following table, where "dim"  means the $k$-dimension of $R/(f)^{m-jsc}$:
\begin{center}
\begin{tabular}{|c|c|c|c|c|c|c|}
\hline 
& \multicolumn{2}{|c|}{$f=x^{3}+y^{4}$} & \multicolumn{2}{|c|}{$f=x^{3}+xy^{3}$} & \multicolumn{2}{|c|}{$f=x^{3}+y^{5}$}\\ 
\hline
m & $(f)^{m-jsc}$ & dim &  $(f)^{m-jsc}$ & dim &  $(f)^{m-jsc}$ & dim    \\
\hline
1& $(x^{2},xy,y^{2})$ & 3 & $(x^{2},xy,y^{2})$ & 3 & $(x^{2},xy,y^{2})$ & 3  \\
\hline
2& $(x^{3},x^{2}y,xy^{2},y^{3})$ & 6 & $(x^{3},x^{2}y,xy^{2},y^{3})$ & 6 & $(x^{3},x^{2}y,xy^{2},y^{3})$ & 6  \\
\hline
3& $(x^{2},xy^{2},y^{4})$ & 6 & $(x^{2},xy^{2},y^{4})$ & 6 & $(x^{2},xy^{2},y^{4})$ & 6 \\
\hline
4& $(x^{3},x^{2}y,xy^{2},y^{3})$ & 6 & $(x^{3},x^{2}y,xy^{3},y^{5})$ & 9 & $(x^{3},x^{2}y,xy^{3},y^{5})$ &  9 \\
\hline
5& $(x^{3},x^{2}y,xy^{2},y^{3})$ & 6 & $(x^{3},x^{2}y,xy^{3},y^{6})$ & 10 & $(x^{3},x^{2}y,xy^{2},y^{3})$ & 6  \\
\hline
6& $(x^{3},x^{2}y,xy^{2},y^{4})$ & 7 & $(x^{3},x^{2}y,xy^{3},y^{7})$ & 11 & $(x^{3},x^{2}y,xy^{2},y^{4})$ & 7  \\
\hline
7& $(x^{3},x^{2}y,xy^{3},y^{4})$ & 8 & $(x^{3},x^{2}y,xy^{3},y^{8})$ & 12 & $(x^{3},x^{2}y,xy^{3},y^{4})$ & 8  \\
\hline
8& & & $(x^{3},x^{2}y^{2},xy^{3},y^{9})$ & 14 & $(x^{3},x^{2}y^{2},xy^{3},y^{5})$ & 10  \\
\hline
9& & & & & $(x^{3},x^{2}y,xy^{3},y^{5})$ & 9  \\
\hline
\end{tabular}
\end{center}
\end{proof}

\section{The Filtration Properties for Jet Closure}
\label{s6}

In this section, we will show that jet closure has a filtration, give some properties and raise some questions about it. We fix $R=k[[x_{1},...,x_{n}]]$ in this section.

\begin{definition}(see \cite{6})
Let $A$ be a commutative ring, a filtration on $A$ is a function $f$, where f takes real values together with $\infty$ and satisfies the following condition: \\
(1) $f(1) \ge 0$; $f(0)= \infty$, \\
(2) $f(x-y) \ge$ Min$(f(x),f(y))$,   \\
(3) $f(xy) \ge f(x)+f(y)$.
\end{definition}

First we show that the jet closure has monotony property.
\begin{proposition}
Let $I \subset R$ be an ideal, then $I^{(m+1)-jc} \subset I^{m-jc}$. 
\end{proposition}
\begin{proof}
For any $f \in I^{(m+1)-jc}$, assume that $f(x_{1}^{(1)}t+...+x_{1}^{(m+1)}t^{m+1},...,x_{n}^{(1)}t+...+x_{n}^{(m+1)}t^{m+1})=f_{1}t+...+f_{m+1}t^{m+1}$(mod $t^{m+2}$) and we have $f_{1},...,f_{m+1} \in (I_{m+1}, \mathbf{m}^{e})$. The structure of $(I_{m+1},\mathbf{m}^{e})$ can be analysed as following. Assume $I$ is generated by $g_{1},...,g_{s}$, and $g_{i}(x_{1}^{(1)}t+...+x_{1}^{(m+1)}t^{m+1},...,x_{n}^{(1)}t+...+x_{n}^{(m+1)}t^{m+1})=g_{i,1}t+...+g_{i,m+1}t^{m+1} (1 \le i \le s)$(mod $t^{m+2}$). In this notation, $(I_{m+1},\mathbf{m}^{e})$ is generated by all $g_{i,j}$ and $\mathbf{m}^{e}$. Noticed that the term $x_{i}^{(m+1)}$ will only contribute to the polynomial $g_{i,m+1}$, so $(I_{m},\mathbf{m}^{e})$ will be generated by $g_{i,j}(1 \le j \le m)$ and $\mathbf{m}^{e}$. Now we have $f_{1},...,f_{m} \in (I_{m+1}, \mathbf{m}^{e})$, and $f_{1},...,f_{m}$ do not contain the variables $x_{i}^{(m+1)}$, so $f_{1},...,f_{m} \in (I_{m}, \mathbf{m}^{e})$. This implies that $f \in I^{m-jc}$.
\end{proof}

\begin{remark}
It is not necessarily true that $I^{(m+1)-jsc} \subset I^{m-jsc}$. We have the following counterexample: take $R=k[x,y]$ and $I=(x^{2}y^{3})$, and we can use Theorem $4.1$ to check that $x^{2}y^{2} \in I^{5-jsc}$ but $x^{2}y^{2} \notin I^{4-jsc}$.
\end{remark}

From \cite{3} we know that in Noetherian case, in general every ideal is arc closed. If the algebra $R/I$ is Artinian, we can get further result:
\begin{proposition}
Suppose $R/I$ is Artinian, then there exists an integer $s$ such that $I$ is $s$-jet closed.  
\end{proposition}
\begin{proof}
From Proposition $1.2$ we know that $I$ is arc closed, so $dim_{k}R/I^{\infty-jc}=dim_{k}R/I < \infty$. By Proposition $6.1$, we have $dim_{k}R/I^{(m+1)-jc} \ge dim_{k}R/I^{m-jc}$, so $\{ dim_{k}R/I^{m-jc} \}$ is a monotone integer sequence with a upper boundary $dim_{k}R/I^{\infty-jc}$, it will be stable after some integer $s$. When $dim_{k}R/I^{s-jc}=dim_{k}R/I^{l-jc}$ for all $l \ge s$, we have $I^{s-jc}=I^{l-jc}$ for all $l \ge s$. By Fact $2.6$, $I=I^{\infty-jc}=I^{s-jc}$.
\end{proof}

By this proposition, for an ideal $I$ such that $R/I$ is Artinian, we can define the minimum integer $s$ such that $I$ is $s$-jet closed.
\begin{definition}
If $R/I$ is Artinian, the jet index $j(I)$ of $I$ is defined to be the minimum integer $s$ such that $I$ is $s$-jet closed. In particular, if $f$ is a polynomial with isolated singularity at $0$, the ideal $J(f):=(\partial_{x_{1}}f,...,\partial_{x_{n}}f)$ satisfies the condition of $R/J(f)$ is Artinian. We say the jet Milnor index of $f$, $j_{\mu}(f)$, to be the jet index of $J(f)$, and we say the jet Tjurina index of $f$, $j_{\tau}(f)$, to be the jet index of $(f,J(f))$.
\end{definition}

By definition, the jet index gives how much of the order of the jet is if we want to give the full information of an ideal. In particular, when the polynomial $f$ is weighted homogeneous, we have $(f,J(f))=J(f)$, so $j_{\mu}(f)=j_{\tau}(f)$. In this situation we can just call them the jet index of $f$, denoted by $j(f)$. If the polynomial $f$ is of the form $x_{1}^{a_{1}+1}+...+x_{n}^{a_{n}+1}$, where $a_{1},...,a_{n}$ are integers, we can give the jet index directly.
\begin{proposition}
The jet index of $f=x_{1}^{a_{1}+1}+...+x_{n}^{a_{n}+1}$ is $a_{1}+...+a_{n}-n$.   
\end{proposition}
\begin{proof}
$J(f)=(x_{1}^{a_{1}},...,x_{n}^{a_{n}})$ is a homogeneous ideal, so by Theorem $3.1$, $$J(f)^{m-jc}=J(f)+\mathbf{m}^{m+1},$$ for any integer $m$. The jet index of $f$ is actually the smallest integer $s$ such that $\mathbf{m}^{s+1} \subset J(f)$, so $j(f)=a_{1}+...+a_{n}-n$.
\end{proof}

\begin{example}
We calculate the jet index of $A_{n},D_{n},E_{6},E_{7},E_{8}$ singularity as examples. 

$A_{n}:f=x^{n+1}+y^{2}$. According to the proposition above, $j(f)=n-1$. 

$E_{6}:f=x^{3}+y^{4}$.  Similarly $j(f)=3$.   

$E_{7}:f=x^{3}+xy^{3}$.  Direct computation with Singular shows that $j(f)=4$ and $J(f)^{m-jc}=J(f)+\mathbf{m}^{m+1}$ for $m=1,2,3$. 

$E_{8}:f=x^{3}+y^{5}$.  Similarly $j(f)=4$. 

$D_{n}:f=x^{n-1}+xy^{2}$.  We will prove that $J(f)^{m-jc}=J(f)+\mathbf{m}^{m+1}$ for $m=1,2,...,n-3$, and $J(f)^{m-jc}=J(f)$ for $m \ge n-2$. 
\end{example}
\begin{proof}
From $f=x^{n-1}+xy^{2}$ we know $J(f)=((n-1)x^{n-2}+y^{2},xy)$. Since $J(f)+\mathbf{m}^{m+1} \subset J(f)^{m-jc}$, $(J(f)+\mathbf{m}^{m+1})^{m-jc} \subset (J(f)^{m-jc})^{m-jc} =J(f)^{m-jc}$, so $(J(f)+\mathbf{m}^{m+1})^{m-jc}=J(f)^{m-jc}$. When $m \le n-3$, $J(f)^{m-jc}=(y^{2},xy)^{m-jc}=(y^{2},xy)+\mathbf{m}^{m+1}=J(f)+\mathbf{m}^{m+1}$ by Proposition $3.6$. When $m=n-2$, suppose $g \notin J(f)$ and $g \in J(f)^{m-jc}$, write $g(x,y)=a_{1}x+...+a_{n-2}x^{n-2}+by$ where $a_{1},...,a_{n-2},b \in k$. Let $x=x_{1}t+...+x_{n-2}t^{n-2},y=y_{1}t+...+y_{n-2}t^{n-2}$ in $g(x,y)$ and we have $a_{1}x_{1}+by_{1} \in (J(f)_{n-2},\mathbf{m}^{e},t^{n-1})$ in $R_{n-2}$. But when module $ (\mathbf{m}^{e},t^{n-1})$, the generators of $J(f)_{n-2}$ have degree at least $2$, so $a_{1}=b=0$. Consider the quotient map:
\begin{flalign*}
q_{1}:R_{n-2}[t]/(J(f)_{n-2},\mathbf{m}^{e},t^{n-1}) \longrightarrow R_{n-2}[t]/(J(f)_{n-2},\mathbf{m}^{e},t^{n-1},y_{1},...y_{n-2}), 
\end{flalign*}
and we have $q_{1} \circ \mu_{n-2}(g)=0$ because $\mu_{n-2}(g)=0$. Through calculating, we can get that $R_{n-2}[t]/(J(f)_{n-2},\mathbf{m}^{e},t^{n-1},y_{1},...y_{n-2})=R_{n-2}[t]/(x_{1}^{n-2},\mathbf{m}^{e},t^{n-1},y_{1},...y_{n-2})$, so $a_{2}=...=a_{n-3}=0$ and $g=a_{n-2}x^{n-2}$. Consider another quotient map:
\begin{flalign*}
q_{2}:R_{n-2}[t]/(J(f)_{n-2},\mathbf{m}^{e},t^{n-1}) \longrightarrow R_{n-2}[t]/(J(f)_{n-2},\mathbf{m}^{e},t^{n-1},y_{2},...y_{n-4},x_{2},...,x_{n-2}) 
\end{flalign*}
and similarly we can get $x_{1}^{n-2} \in (x_{1}^{n-2}+2y_{1}y_{n-3},y_{1}^{2},x_{1}y_{1},x_{1}y_{n-3})$, so $y_{1}y_{n-3} \in (y_{1}^{2},x_{1}y_{1},x_{1}y_{n-3})$, this is a contradiction.
\end{proof}

It is natural to compare the jet Milnor index $j_{\mu}(f)$ with index of nilpotency of the maximal ideal of the Artinian ring $R/J(f)$, which we denote by $N(J(f))$. The latter satisfies $\mathbf{m}^{N(J(f))} \subset J(f)$ and $\mathbf{m}^{N(J(f))-1} \nsubseteq J(f)$. Since $\mathbf{m}^{j_{\mu}(f)+1} \subset J(f)^{j_{\mu}(f)-jc}=J(f)$, we have  $j_{\mu}(f)+1 \ge N(J(f))$. This inequality can be strict, so the concept of jet Milnor index is different from that of nilpotency of maximal ideal of the Artinian ring. The example is as follows, which also tells us the equation $J(f)^{m-jc}=J(f)+\mathbf{m}^{m+1}$ does not always hold.  

\begin{example}
Take $f=x^{6}+xy^{6}+x^{2}y^{5}$, and $J(f)=(6x^{5}+y^{6}+2xy^{5},6xy^{5}+5x^{2}y^{4})$. Direct computation shows the common zero locus of $f$ and $J(f)$ is $0$, so $f$ has an isolated singularity at $0$. One can use Singular to check that $\mathbf{m}^{10} \nsubseteq J(f)$ but $\mathbf{m}^{11} \subset J(f)$, so $N(J(f))=11$. Next we will show that $j_{\mu}(f)>10$, i.e. $J(f) \neq J(f)^{10-jc}$. In fact, one can use definition to check directly that $x^{8} \in J(f)^{10-jc}$, but with Singular one can find that $x^{8} \notin J(f)+\mathbf{m}^{11}$. This implies $j_{\mu}(f)>N(J(f))-1$ and $J(f)^{10-jc} \neq J(f)+\mathbf{m}^{11}$.
\end{example}

Similarly we have another inequality $j_{\tau}(f)+1 \ge N(f,J(f))$. However, through a large amount of calculations we can not find an example such that the equality does not hold, so we propose the following conjecture:
\begin{conjecture}
Suppose $f \in k[[x_{1},...,x_{n}]]$ is a polynomial which defines an isolated singularity at $0$, then we have $j_{\tau}(f)+1=N(f,J(f))$.
\end{conjecture}

We can show that the gap between nearby ideals in $\{I^{m-jc}\}_{m=0}^{\infty}$ is within the maximal ideal $\mathbf{m}$. Before that, we prepare a basic lemma about jet ideal $I_{m}$.
\begin{lemma}
If we embed $R_{m}$ into $R_{m+1}$ canonically, we have $I_{m} \subset I_{m+1}$, and thus $\sqrt{I_{m}} \subset \sqrt{I_{m+1}}$.  
\end{lemma}
\begin{proof}
Since $(R/I)_{m}=R_{m}/I_{m}$, we have an injective map $R_{m}/I_{m} \longrightarrow R_{m+1}/I_{m+1}$, so $I_{m} \subset I_{m+1}$.
\end{proof}

We can use this lemma to prove the following proposition.

\begin{proposition}
For any ideal $I \subset R$, we have $\mathbf{m}I^{m-jc} \subset I^{(m+1)-jc}$.
\end{proposition}  
\begin{proof}
It suffices to show for any $f \in I^{m-jc}$, $x_{1}f \in I^{(m+1)-jc}$. Suppose $$f(x_{1}^{(1)}t+...+x_{1}^{(m)}t^{m},...,x_{n}^{(1)}t+...+x_{n}^{(m)}t^{m})=f_{1}t+...+f_{m}t^{m}+...+f_{l}t^{l}$$ for some integer $l$ and we have $f_{1},...,f_{m} \in I_{m}$. Now we calculate $\mu_{m+1}(x_{1}f)=(x_{1}^{(1)}t+...+x_{1}^{(m+1)}t^{m+1})f(x_{1}^{(1)}t+...+x_{1}^{(m+1)}t^{m+1},...,x_{n}^{(1)}t+...+x_{n}^{(m+1)}t^{m+1})$ in $R_{m+1}[t]/(t^{m+2},\mathbf{m}^{e},I_{m+1})$. Since $f(0)=0$,
\begin{flalign*}
&(x_{1}^{(1)}t+...+x_{1}^{(m+1)}t^{m+1})f(x_{1}^{(1)}t+...+x_{1}^{(m+1)}t^{m+1},...,x_{n}^{(1)}t+...+x_{n}^{(m+1)}t^{m+1}) \\
=&(x_{1}^{(1)}t+...+x_{1}^{(m)}t^{m})f(x_{1}^{(1)}t+...+x_{1}^{(m)}t^{m},...,x_{n}^{(1)}t+...+x_{n}^{(m)}t^{m})  \\
=&(x_{1}^{(1)}t+...+x_{1}^{(m)}t^{m})(f_{1}t+...+f_{m}t^{m}) (mod (I_{m+1},\mathbf{m}^{e},t^{m+2})).
\end{flalign*}

Because $f_{1},...,f_{m} \in I_{m} \subset I_{m+1}$ by Lemma $6.9$, each coefficient of $t^{i}(1 \le i \le m+1)$ is in $(I_{m+1},\mathbf{m}^{e},t^{m+1})$. This shows $\mu_{m+1}(x_{1}f)=0$ in $R_{m+1}[t]/(t^{m+1},\mathbf{m}^{e},I_{m+1})$, and the proposition follows.
\end{proof}

In order to study $\{I^{m-jc}\}_{m=0}^{\infty}$, we define a map:
\begin{flalign*}
f_{I}: R/I \longrightarrow \mathbb{Z}_{\ge 0}.
\end{flalign*}
For any $g \in R/I$, if $g \in I^{m-jc}/I$ and $g \notin I^{(m+1)-jc}/I$, take $f_{I}(g)=m+1$. If $g \in I^{\infty-jc}/I$, take $f_{I}(g)=\infty $. If we set $I^{0-jc}=I+\mathbf{m}=\mathbf{m}$, and if $g \notin \mathbf{m}$, take $f_{I}(g)=0$, we get a map satisfying the following properties:
\begin{proposition}
(1)$f_{I}(1)=0, f_{I}(0)=\infty$;   

(2)$f_{I}(x+y) \ge min \{ f_{I}(x), f_{I}(y) \}$;  

(3)$f_{I}(xy) \ge f_{I}(x)+f_{I}(y)$.
\end{proposition}
\begin{proof}
(1) Directly from the definition. 

(2)If $min \{ f_{I}(x), f_{I}(y) \}=\infty$, $x,y \in I^{\infty-jc}$, so $x+y \in I^{\infty-jc}$ and the inequality holds. If $min \{ f_{I}(x), f_{I}(y) \}=m$, then $x,y \in I^{(m-1)-jc}$, so $x+y \in I^{(m-1)-jc}$ and the inequality holds. 

(3)It suffices to show for any $p,q \in \mathbb{Z}_{\ge 0}$, $I^{p-jc}I^{q-jc} \subset I^{(p+q+1)-jc}$. For any $f \in I^{p-jc}$, $g \in I^{q-jc}$, assume $f(x_{1}^{(1)}t+...+x_{1}^{(p+q+1)}t^{p+q+1},...,x_{n}^{(1)}t+...+x_{n}^{(p+q+1)}t^{p+q+1})=f_{1}t+...+f_{p}t^{p}+...$, and $g(x_{1}^{(1)}t+...+x_{1}^{(p+q+1)}t^{p+q+1},...,x_{n}^{(1)}t+...+x_{n}^{(p+q+1)}t^{p+q+1})=g_{1}t+...+g_{q}t^{q}+...$, then $f_{1},...,f_{p} \in I_{p}$ and $g_{1},...,g_{q} \in I_{q}$. Now we calculate $f*g$ as follows:
\begin{flalign*}
&f*g(x_{1}^{(1)}t+...+x_{1}^{(p+q+1)}t^{p+q+1},...,x_{n}^{(1)}t+...+x_{n}^{(p+q+1)}t^{p+q+1}) \\
=&(f_{1}t+...+f_{p}t^{p}+...)(g_{1}t+...+g_{q}t^{q}+...) \\
=&h_{1}t+...+h_{p+q+1}t^{p+q+1}.
\end{flalign*}

Notice that for any $1 \le i \le p+q+1$, each monomial of $h_{i}$ has a factor which has the form $f_{j}(1 \le j \le p)$ or $g_{l}(1 \le l \le q)$, so every $h_{i}$ is in $I_{p+q+1}$, this gives the assertion.
\end{proof}

This proposition tells us $f_{I}$ is a filtration on $R/I$. The author of \cite{6} also shows some general properties about filtration, one of them, for example, we can explore is if the filtration is homogeneous.

\begin{definition}
A filtration $f: R/I \longrightarrow \mathbb{Z}_{\ge 0}$ is called homogeneous if for all $s \in \mathbb{N},x \in R$, we have $f(x^{s})=sf(x)$.
\end{definition}

Since the jet closure of homogeneous ideal and monomial ideal is easy to compute, we can get directly if their filtration is homogeneous.

\begin{proposition}
If the ideal $I$ is homogeneous or it is a monomial ideal, the filtration $f_{I}$ is homogeneous if and only if $I$ is a radical ideal.  
\end{proposition}

\begin{proof}
On one hand, suppose $I$ is a radical ideal. By Theorem A, when $I$ is a homogeneous ideal or monomial ideal, $I^{m-jc}=I+\mathbf{m}^{m+1}$. When $x \in I$, $f_{I}(x)=\infty$ and $f_{I}(x^{s})=\infty$, so the equality holds. Suppose $f_{I}(x)=r$, we have $x \in I+\mathbf{m}^{r}/I$ and $x \notin I+\mathbf{m}^{r+1}/I$, so we can take $x$ as a homogeneous polynomial with degree $m$. If $f(x^{s}) \ge rs+1$, we have $x^{s} \in I+\mathbf{m}^{rs+1}/I$. Assume now $x^{s}=a+b$, where $a \in I$ and $b \in \mathbf{m}^{rs+1}$, so $x^{s}-b \in I$. Since $I$ is a homogeneous ideal or monomial ideal, and $x^{s}$ is a homogeneous polynomial of degree $rs$, we have $x^{s} \in I$. Due to $I$ is a radical ideal, $x \in I$, this is a contradiction. 

On the other hand, if $f_{I}$ is homogeneous. For any $x \in \sqrt{I}$, there exists an integer $s$ such that $x^{s} \in I$, so $\infty=f_{I}(x^{s})=sf(x)$. This implies that $f_{I}(x)=\infty$, so $x \in I$. 
\end{proof}

When the ideal $I$ does not have such properties, how to characterize the  filtration ${f}_{I}$ is a good question. Here I present the following question.

\begin{question}
Generally when is the filtraion $f_{I}$ is homogeneous? 
\end{question}

Lemma $2.11$ in \cite{6} tells us that for any filtration $f$, the limit of $f(x^{n})/n$ always exists as $n$ increase to $\infty$, and if we write this limit $\mathbf{f}(x)$, we can get a homogeneous filtration $\mathbf{f}$ with $\mathbf{f} \ge f$. Thus for the filtration $f_{I}$, we can get a new homogeneous filtration $\mathbf{f}_{I}$. We can compute it when $I$ is a homogeneous ideal or monomial ideal with the following proposition.

\begin{proposition}
For any homogeneous ideal or monomial ideal $I$, the filtration $\mathbf{f}_{I}$ is given by $\{ \sqrt{I}+\mathbf{m}^{m} \}_{m=1}^{\infty}$. 
\end{proposition}
\begin{proof}
For any $x \in R/I$, if $x \in \sqrt{I}$, then there exists an integer $s$ such that $x^{s} \in I$, so $\mathbf{f}(x)=\infty$. If $x \notin \mathbf{m}$, then $x^{r} \notin \mathbf{m}$ for any integer $r$, so $\mathbf{f}(x)=0$. Now suppose $x \in \sqrt{I}+\mathbf{m}^{r}/I$ and $x \notin \sqrt{I}+\mathbf{m}^{r+1}/I$. We assume that $x=a+b$, where $a \in \sqrt{I}$, $b \in \mathbf{m}^{r}$ and $a^{s} \in I$. We have $x^{m}=(a+b)^{m} \in I+\mathbf{m}^{(m-s+1)r}$, so $f(x^{n}) \ge (m-s+1)r$. This tells us that $\mathbf{f}(x) \ge r$. If $\mathbf{f}(x) >r$, then there exists an integer $u$ such that $f(x^{u}) \ge ur+1$, so $x^{u} \in I+\mathbf{m}^{ur+1}/I \subset \sqrt{I}+\mathbf{m}^{ur+1}/I$. Now consider the expansion of $x^{u}=(a+b)^{u}$: since $a \in \sqrt{I}$, we have $x^{u}-b^{u} \in \sqrt{I}$. Notice that $b^{u} \in \mathbf{m}^{ur}$ and $\sqrt{I}$ is a homogeneous ideal or monomial ideal because $I$ is a homogeneous ideal or monomial ideal, so $x^{n} \in \sqrt{I}$, this is a contradiction. This tells us $\mathbf{f}(x)=r$, which is the required conclusion.
\end{proof}

\section{Appendix: Codes for Calculation}
\label{s7}

Here we list two codes in Singular for computing jet closure and jet support closure. We use Corollary $2.1$ and Proposition $2.5 (3)$ to calculate. First we give the code for computing jet closure:

\begin{algorithm}
(Calculating jet closure) \\
proc jetclosure(ideal I,int m)         \\
 \{                                      \\
  def R=basering;                       \\
  int n=nvars(R);                       \\
  int s=size(I);                        \\
  ring Rhelp=0,(x(1..n)(0..m),t),dp;    \\
  poly f(1..n)=x(1..n)(0);              \\
  for(int i=1;i $\le$ n;i++)                 \\
  \{                                     \\
   for(int j=1;j $\le$m;j++)                \\
   \{                                    \\
    f(i)=f(i)+x(i)(j)*$t^{j}$;              \\
   \};                                   \\
  \};                                    \\
 map F=R, f(1..n);                      \\
 ideal J;                               \\
 for(int i=1;i$\le$s;i++)                  \\
 \{                                      \\
  poly g=F(I)[i];                       \\
  for(int j=0;j$\le$m;j++)                 \\
  \{                                     \\
   ideal K=g,t;                         \\
   ideal L=eliminate(K,t);              \\
   int a=(i-1)*(m+1)+j+1;               \\
   J[a]=L[1];                           \\
   g=g-L[1];                            \\
   g=g/t;                               \\
  \};                                    \\
 \};                                     \\
 ideal M=J,x(1..n)(0),$t^{m+1}$;          \\
 setring R;                             \\
 ideal N=preimage(Rhelp,F,M);           \\
return (N);                             \\
 \};
\end{algorithm}

The code for computing jet support closure is as follows.
\begin{algorithm}
(Calculating jet support closure) \\
LIB"primdec.lib";              \\
proc jetsuppclosure(ideal I,int m)         \\
 \{                                      \\
  def R=basering;                       \\
  int n=nvars(R);                       \\
  int s=size(I);                        \\
  ring Rhelp=0,(x(1..n)(0..m),t),dp;    \\
  poly f(1..n)=x(1..n)(0);              \\
  for(int i=1;i $\le$ n;i++)                 \\
  \{                                     \\
   for(int j=1;j $\le$m;j++)                \\
   \{                                    \\
    f(i)=f(i)+x(i)(j)*$t^{j}$;              \\
   \};                                   \\
  \};                                    \\
 map F=R, f(1..n);                      \\
 ideal J;                               \\
 for(int i=1;i$\le$s;i++)                  \\
 \{                                      \\
  poly g=F(I)[i];                       \\
  for(int j=0;j$\le$m;j++)                 \\
  \{                                     \\
   ideal K=g,t;                         \\
   ideal L=eliminate(K,t);              \\
   int a=(i-1)*(m+1)+j+1;               \\
   J[a]=L[1];                           \\
   g=g-L[1];                            \\
   g=g/t;                               \\
  \};                                    \\
 \};                                    \\
 J=J,x(1..n)(0);                                     \\
 ideal M=radical(J),$t^{m+1}$;          \\
 setring R;                             \\
 ideal N=preimage(Rhelp,F,M);           \\
return (N);                             \\
 \};
\end{algorithm}

\end{document}